\documentclass{amsart}
\usepackage{geometry,amssymb,color,hyperref}
\usepackage{graphicx} 

\theoremstyle{plain}

\theoremstyle{definition}

\newcommand{\IZ}{\mathbb Z}

\theoremstyle{remark}

\newtheorem{example}{Example}[section]

\title{New Steiner systems $S(2,6,v)$ with block length 6}
\author{Taras Banakh, Ivan Hetman, Alex Ravsky}
\date{July 2025}

\begin{document}

    \begin{abstract} In this paper various Steiner systems $S(2,k,v)$ for $k = 6$ are collected and enumerated for specific constructions. In particular, two earlier unknown types of $1$-rotational designs are found for the groups $SL(2,5)$ and $((\mathbb Z_3 \times \mathbb Z_3) \rtimes \mathbb Z_3) \times \mathbb Z_5$. Also new Steiner systems $S(2,6,66)$, $S(2,6,96), S(2,6,106), S(2,6,111)$ are listed.
    \end{abstract}

    \maketitle

    \section{Introduction}

    This paper describes some Steiner systems $S(2,6,v)$ obtained by various algorithms. The paper \cite{Het} enumerates $1$-rotational Steiner systems $S(2,k,v)$ for $k\in\{3,4,5\}$, whereas \cite{Het1} enumerates 1-rotational unitals, which are Steiner systems $S(2,6,126)$. The next obvious step would be to try finding other $1$-rotational for $k\ge 6$. While for the case $k \le 5$ our  algorithm still has capacity of generating new enumerations, for $k \ge 6$ the algorithm used in \cite{Het} starts struggling with computational complexity. Fortunately, a mixed approach uniting strong sides of a generalized algorithm and the algorithm for commutative groups was elaborated, and this paper contains  results on $1$-rotational designs for $k = 6$ found by this ``mixed'' algorithm.
    Surprisingly, in addition to designs found by W.~H.~Mills \cite{Mills} for the group {\tt C5\,x\,(C9\,:\,C3)}, two new earlier unknown types of $1$-rotational designs were found and enumerated for the groups {\tt SL(2,5)} and  {\tt C5\,x\,((C3\,x\,C3)\,:\,C3)}.

    By investigating different acting groups of order $20$ with 5 orbits of size $20, 20, 20, 5, 1$, it was possible to find one new design $S(2,6,66)$ with acting group $\mathbb Z_{2} \times \mathbb Z_{10}$.

    In addition, one more algorithm was implemented that covers the case of two orbits of size $n$ with acting group $\mathbb Z_{n}$. This allowed us to enumerate designs obtained by Mills in \cite{Mills1} \cite{Mills2} and in addition to find four new designs $S(2,6,111)$ with acting group $\mathbb Z_{55}$.

    The structure analysis showed that there exists a unique design generated by the group $\mathbb Z_{48}$. Nevertheless, there is an interesting  structure introduced by Denniston \cite{Denniston} who generated a Steiner system $S(2,6,66)$ using the acting group $\mathbb Z_{13} \rtimes \mathbb Z_{3}$ with four orbits of size $39, 13, 13, 1$. We have studied a similar structure for the acting group $\mathbb Z_{19} \rtimes \mathbb Z_{3}$ on four orbits of size $57, 19, 19, 1$ and found that it generates exactly 133 Steiner designs $S(2,6,96)$.

    \section{1-rotational designs with $k = 6$}

    In this section we present some Steiner systems generated by the left action of a group $G$ on $G$ with added fixed point. Such designs are called {\em $1$-rotational Steiner systems}, see VI.16.6 in \cite{HoCD}. The following table summarizes the results of our calculations.
    \smallskip

    \begin{center}
        \begin{tabular}{|c c c c|}
            \hline
            \phantom{$|^{|^|}$} {GAP ID}\phantom{$|^{|^|}$} & {group structure} & {number of designs} & {comments} \\
            \hline
            {\tt SmallGroup(120,5)} & {\tt SL(2,5)} & 24 & - \\
            \hline
            {\tt SmallGroup(125,1)} & $\mathbb Z_{125}$ & 8 & \cite{Het1} \\
            {\tt SmallGroup(125,2)} & $\mathbb Z_5 \times \mathbb Z_{25}$ & 32 & \cite{Het1} \\
            {\tt SmallGroup(125,3)} & $(\mathbb Z_5 \times \mathbb Z_5) \rtimes \mathbb Z_5$ & 20 & \cite{Het1} \\
            {\tt SmallGroup(125,4)} & $\mathbb Z_{25} \rtimes \mathbb Z_5$ & 29 & \cite{Het1} \\
            {\tt SmallGroup(125,5)} & $\mathbb Z_5 \times \mathbb Z_5 \times \mathbb Z_5$ & 8 & \cite{Het1} \\
            \hline
            {\tt SmallGroup(135,3)} & {\tt C5\,x\,((C3\,x\,C3)\,:\,C3)} & 11 & - \\
            {\tt SmallGroup(135,4)} & {\tt C5\,x\,(C9\,:\,C3)} & 20 & \cite{Mills} \\
            \hline
            {\tt SmallGroup(150,5)} & {\tt (C5\,x\,C5)\,:\,S3} & ? & - \\
            {\tt SmallGroup(150,7)} & {\tt C2\,x\,((C5\,x\,C5)\,:\,C3)} & $\ge$ 4 & - \\
            {\tt SmallGroup(150,11)} & {\tt C5\,x\,D30} & ? & - \\
            \hline
            {\tt SmallGroup(155,1)} & {\tt C31\,:\,C5} & $\ge$ 1 & - \\
            {\tt SmallGroup(155,2)} & $\mathbb Z_{155}$ & $\ge$ 16 & - \\
            \hline
        \end{tabular}
    \end{center}

    The $1$-rotational designs are presented in the form of a list of basic blocks on which the group acts from the left, thus generating the whole design. Each design is preceded by its {\em fingerprint}, described in \cite{Het} \cite{Het1}. The fingerprint is calculated in  $O(k^2 v^3)$ time and effectively distinguishes non-isomorphic designs in the sense that two designs with different fingerprints are certainly non-isomorphic. If for two designs their fingerprinets coincide then a more elaborated algorithms (like {\tt Nauty}) should be applied for determining the (non)isomorphness of these designs.

    For calculation purposes it is necessary to transform a group from the GAP representation to its Cayley table. The simplest way is to use LOOPS \cite{LOOPS} package for GAP applying the procedure {\tt CayleyTable(IntoLoop(SmallGroup(126,1)));} to obtain a two-dimensional array and then convert it to 0-based Cayley table by subtracting $1$. Then the generated Cayley table is used for producing all blocks from the basic blocks.

    \begin{example} There are no 1-rotational designs for $v\in\{31, 36, 46, 51, 61, 66, 76, 81, 91, 96, 106, 111,\allowbreak 141, 166\}$.
    \end{example}

    \begin{example} {\tt SmallGroup(120,5) = SL(2,5)}
        \begin{enumerate}
            \item \{1=24000, 2=601200, 3=2665440, 4=3388560\} \newline [[0, 1, 2, 3, 4, $\infty$], [0, 5, 18, 49, 59, 105], [0, 6, 44, 55, 68, 75], [0, 7, 42, 45, 60, 111], \newline [0, 11, 16, 30, 41, 61], [0, 21, 52, 69, 89, 90], [0, 24, 50, 67, 78, 104]]
            \item \{1=24960, 2=603360, 3=2744640, 4=3306240\} \newline [[0, 1, 2, 3, 4, $\infty$], [0, 5, 19, 34, 42, 106], [0, 6, 26, 61, 89, 93], [0, 7, 12, 38, 57, 113], \newline [0, 10, 48, 49, 58, 84], [0, 13, 28, 56, 90, 103], [0, 16, 67, 72, 79, 104]]
            \item \{1=30720, 2=603360, 3=2761920, 4=3283200\} \newline [[0, 1, 2, 3, 4, $\infty$], [0, 5, 15, 23, 24, 77], [0, 7, 46, 92, 98, 117], [0, 9, 45, 68, 73, 112], \newline [0, 10, 52, 57, 89, 104], [0, 11, 16, 30, 41, 61], [0, 42, 43, 70, 94, 96]]
            \item \{1=32640, 2=576720, 3=2726880, 4=3342960\} \newline [[0, 1, 2, 3, 4, $\infty$], [0, 5, 11, 27, 30, 55], [0, 6, 65, 71, 89, 111], [0, 7, 29, 54, 63, 91], \newline [0, 8, 10, 48, 68, 98], [0, 16, 67, 72, 79, 104], [0, 23, 24, 66, 92, 115]]
            \item \{1=32640, 2=594720, 3=2731200, 4=3320640\} \newline [[0, 1, 2, 3, 4, $\infty$], [0, 5, 32, 82, 99, 104], [0, 6, 12, 15, 59, 61], [0, 7, 46, 55, 84, 119], \newline [0, 8, 75, 109, 111, 114], [0, 16, 23, 28, 87, 90], [0, 42, 43, 70, 94, 96]]
            \item \{1=33600, 2=600480, 3=2759040, 4=3286080\} \newline [[0, 1, 2, 3, 4, $\infty$], [0, 5, 15, 23, 24, 77], [0, 7, 46, 92, 98, 117], [0, 9, 34, 59, 81, 116], \newline [0, 10, 52, 57, 89, 104], [0, 11, 16, 30, 41, 61], [0, 42, 43, 70, 94, 96]]
            \item \{1=33600, 2=626400, 3=2729280, 4=3289920\} \newline [[0, 1, 2, 3, 4, $\infty$], [0, 5, 32, 82, 99, 104], [0, 6, 23, 79, 80, 94], [0, 7, 12, 46, 63, 66], \newline [0, 8, 70, 105, 111, 119], [0, 11, 15, 30, 34, 110], [0, 16, 77, 102, 108, 114]]
            \item \{1=34560, 2=576000, 3=2732160, 4=3336480\} \newline [[0, 1, 2, 3, 4, $\infty$], [0, 5, 29, 37, 59, 101], [0, 6, 41, 80, 110, 111], [0, 8, 24, 50, 81, 86], \newline [0, 10, 48, 49, 58, 84], [0, 12, 16, 42, 65, 70], [0, 22, 28, 68, 90, 91]]
            \item \{1=35520, 2=592560, 3=2701920, 4=3349200\} \newline [[0, 1, 2, 3, 4, $\infty$], [0, 5, 14, 69, 83, 116], [0, 6, 55, 92, 94, 101], [0, 7, 62, 63, 96, 115], \newline [0, 10, 11, 17, 30, 66], [0, 16, 23, 28, 87, 90], [0, 18, 49, 56, 58, 117]]
            \item \{1=36480, 2=576720, 3=2747040, 4=3318960\} \newline [[0, 1, 2, 3, 4, $\infty$], [0, 5, 11, 27, 30, 55], [0, 6, 61, 101, 104, 109], [0, 7, 47, 65, 71, 108], \newline [0, 9, 40, 59, 110, 113], [0, 15, 52, 58, 89, 116], [0, 16, 23, 28, 87, 90]]
            \item \{1=41280, 2=606240, 3=2705280, 4=3326400\} \newline [[0, 1, 2, 3, 4, $\infty$], [0, 5, 14, 59, 67, 74], [0, 6, 55, 78, 101, 106], [0, 9, 45, 56, 94, 105], \newline [0, 10, 48, 49, 58, 84], [0, 11, 30, 42, 69, 116], [0, 16, 32, 62, 109, 113]]
            \item \{1=42240, 2=603360, 3=2699520, 4=3334080\} \newline [[0, 1, 2, 3, 4, $\infty$], [0, 5, 44, 95, 113, 116], [0, 6, 27, 31, 51, 78], [0, 11, 12, 30, 73, 96], \newline [0, 13, 19, 23, 60, 83], [0, 14, 49, 55, 58, 102], [0, 16, 67, 72, 79, 104]]
            \item \{1=47040, 2=596160, 3=2740800, 4=3295200\} \newline [[0, 1, 2, 3, 4, $\infty$], [0, 5, 32, 82, 99, 104], [0, 6, 39, 40, 85, 92], [0, 7, 33, 78, 111, 112], \newline [0, 8, 57, 65, 80, 110], [0, 11, 16, 30, 41, 61], [0, 28, 53, 90, 94, 107]]
            \item \{1=49920, 2=597600, 3=2729280, 4=3302400\} \newline [[0, 1, 2, 3, 4, $\infty$], [0, 5, 11, 27, 30, 55], [0, 6, 23, 81, 95, 119], [0, 7, 43, 82, 99, 112], \newline [0, 8, 44, 46, 51, 80], [0, 9, 61, 77, 93, 116], [0, 16, 49, 54, 58, 83]]
            \item \{1=49920, 2=601200, 3=2696160, 4=3331920\} \newline [[0, 1, 2, 3, 4, $\infty$], [0, 5, 34, 38, 50, 76], [0, 6, 12, 13, 57, 84], [0, 9, 31, 81, 86, 111], \newline [0, 10, 11, 17, 30, 66], [0, 16, 32, 62, 109, 113], [0, 20, 49, 58, 63, 73]]
            \item \{0=1200, 1=25920, 2=570960, 3=2713440, 4=3367680\} \newline [[0, 1, 2, 3, 4, $\infty$], [0, 5, 11, 27, 30, 55], [0, 6, 12, 43, 44, 63], [0, 7, 18, 45, 96, 112], \newline [0, 8, 14, 90, 102, 104], [0, 10, 48, 49, 58, 84], [0, 16, 57, 82, 99, 119]]
            \item \{0=1200, 1=28800, 2=602640, 3=2696160, 4=3350400\} \newline [[0, 1, 2, 3, 4, $\infty$], [0, 5, 57, 77, 100, 112], [0, 6, 63, 67, 107, 113], [0, 7, 29, 33, 34, 102], \newline [0, 10, 11, 17, 30, 66], [0, 16, 49, 54, 58, 83], [0, 22, 28, 68, 90, 91]]
            \item \{0=1200, 1=32640, 2=582480, 3=2734560, 4=3328320\} \newline [[0, 1, 2, 3, 4, $\infty$], [0, 5, 38, 48, 93, 101], [0, 7, 51, 57, 95, 116], [0, 9, 13, 46, 47, 68], \newline [0, 10, 11, 17, 30, 66], [0, 12, 52, 53, 89, 110], [0, 16, 49, 54, 58, 83]]
            \item \{0=1200, 1=32640, 2=603360, 3=2702400, 4=3339600\} \newline [[0, 1, 2, 3, 4, $\infty$], [0, 5, 19, 33, 65, 102], [0, 6, 46, 54, 100, 112], [0, 10, 48, 49, 58, 84], \newline [0, 11, 30, 63, 82, 88], [0, 14, 17, 66, 81, 104], [0, 16, 23, 28, 87, 90]]
            \item \{0=1200, 1=35520, 2=591120, 3=2693280, 4=3358080\} \newline [[0, 1, 2, 3, 4, $\infty$], [0, 5, 38, 42, 80, 96], [0, 6, 14, 34, 39, 110], [0, 10, 48, 49, 58, 84], \newline [0, 11, 30, 63, 82, 88], [0, 12, 44, 61, 73, 83], [0, 16, 67, 72, 79, 104]]
            \item \{0=1200, 1=39360, 2=571680, 3=2716800, 4=3350160\} \newline [[0, 1, 2, 3, 4, $\infty$], [0, 5, 11, 27, 30, 55], [0, 6, 74, 83, 105, 107], [0, 7, 28, 61, 67, 90], \newline [0, 8, 14, 88, 96, 115], [0, 10, 29, 44, 63, 80], [0, 16, 32, 62, 109, 113]]
            \item \{0=1200, 1=41280, 2=624240, 3=2711520, 4=3300960\} \newline [[0, 1, 2, 3, 4, $\infty$], [0, 5, 11, 12, 19, 116], [0, 6, 46, 101, 110, 118], [0, 7, 24, 48, 59, 79], \newline [0, 10, 52, 57, 89, 104], [0, 16, 49, 54, 58, 83], [0, 23, 42, 51, 70, 77]]
            \item \{0=2400, 1=26880, 2=558720, 3=2771520, 4=3319680\} \newline [[0, 1, 2, 3, 4, $\infty$], [0, 5, 11, 27, 30, 55], [0, 6, 51, 53, 68, 92], [0, 7, 80, 87, 96, 98], \newline [0, 8, 10, 48, 75, 104], [0, 15, 52, 58, 89, 116], [0, 16, 24, 29, 95, 115]]
            \item \{0=2400, 1=35520, 2=658800, 3=2687520, 4=3294960\} \newline [[0, 1, 2, 3, 4, $\infty$], [0, 5, 17, 20, 74, 87], [0, 6, 80, 86, 90, 106], [0, 8, 32, 53, 63, 94], \newline [0, 11, 12, 30, 73, 96], [0, 16, 24, 29, 95, 115], [0, 21, 62, 100, 109, 114]]
        \end{enumerate}
    \end{example}

    \begin{example} {\tt SmallGroup(135,3) = C5\,x\,((C3\,x\,C3)\,:\,C3)}
        \begin{enumerate}
            \item \{1=23760, 2=666630, 3=3590460, 4=5266350\}\newline [[0, 1, 13, 17, 57, 101], [0, 3, 12, 29, 52, $\infty$], [0, 4, 70, 97, 122, 132], [0, 11, 16, 33, 48, 55],\newline [0, 15, 39, 95, 96, 107], [0, 19, 44, 60, 89, 91]]
            \item \{1=25920, 2=598590, 3=3601260, 4=5321430\}\newline [[0, 1, 6, 13, 97, 104], [0, 3, 12, 29, 52, $\infty$], [0, 4, 46, 47, 94, 112], [0, 7, 79, 80, 85, 125],\newline  [0, 8, 25, 32, 37, 65], [0, 11, 54, 57, 81, 95]]
            \item \{1=28080, 2=560520, 3=3579120, 4=5379480\}\newline [[0, 2, 4, 21, 108, 120], [0, 3, 12, 29, 52, $\infty$], [0, 6, 51, 76, 125, 133], [0, 7, 13, 66, 91, 116],\newline [0, 8, 37, 75, 114, 132], [0, 10, 15, 78, 81, 96]]
            \item \{1=33480, 2=680400, 3=3660120, 4=5173200\}\newline [[0, 1, 5, 47, 64, 81], [0, 2, 4, 55, 76, 134], [0, 3, 12, 29, 52, $\infty$], [0, 6, 48, 51, 57, 110],\newline [0, 7, 18, 26, 61, 101], [0, 13, 90, 109, 116, 129]]
            \item \{1=37800, 2=679590, 3=3598020, 4=5231790\}\newline [[0, 1, 10, 27, 93, 94], [0, 2, 9, 32, 34, 37], [0, 3, 12, 29, 52, $\infty$], [0, 6, 44, 46, 103, 130],\newline  [0, 11, 24, 64, 96, 111], [0, 15, 100, 107, 113, 116]]
            \item \{1=41040, 2=640710, 3=3603420, 4=5262030\}\newline [[0, 1, 2, 13, 32, 73], [0, 3, 12, 29, 52, $\infty$], [0, 4, 44, 55, 123, 124], [0, 6, 38, 46, 53, 81],\newline  [0, 8, 16, 30, 51, 83], [0, 11, 48, 89, 93, 96]]
            \item \{1=41040, 2=665820, 3=3628800, 4=5211540\}\newline [[0, 1, 2, 50, 87, 128], [0, 3, 12, 29, 52, $\infty$], [0, 4, 41, 97, 111, 132], [0, 7, 63, 74, 98, 127],\newline [0, 11, 48, 89, 93, 96], [0, 13, 26, 32, 91, 115]]
            \item \{0=1350, 1=33480, 2=664200, 3=3519720, 4=5328450\}\newline [[0, 1, 4, 26, 35, 58], [0, 3, 12, 29, 52, $\infty$], [0, 7, 78, 90, 104, 132], [0, 10, 13, 63, 81, 110],\newline [0, 11, 48, 62, 85, 105], [0, 15, 46, 79, 121, 129]]
            \item \{0=1350, 1=33480, 2=720090, 3=3650940, 4=5141340\}\newline [[0, 1, 2, 17, 76, 126], [0, 3, 12, 29, 52, $\infty$], [0, 4, 40, 60, 97, 98], [0, 6, 45, 79, 86, 128],\newline  [0, 11, 48, 89, 93, 96], [0, 16, 20, 71, 118, 119]]
            \item \{0=1350, 1=37800, 2=660960, 3=3555360, 4=5291730\}\newline [[0, 1, 2, 44, 114, 118], [0, 3, 12, 29, 52, $\infty$], [0, 4, 19, 46, 57, 102], [0, 6, 10, 22, 62, 86],\newline  [0, 11, 43, 48, 112, 115], [0, 17, 70, 75, 113, 120]]
            \item \{0=1350, 1=45360, 2=627750, 3=3644460, 4=5228280\}\newline [[0, 1, 4, 32, 49, 77], [0, 3, 12, 29, 52, $\infty$], [0, 6, 24, 31, 116, 130], [0, 7, 39, 61, 72, 100],\newline [0, 8, 41, 88, 109, 115], [0, 11, 48, 62, 85, 105]]
        \end{enumerate}
    \end{example}

    \begin{example} {\tt SmallGroup(135,4) = C5\,x\,(C9\,:\,C3)}\newline The design (17) was found by Mills \cite{Mills}.
    \begin{enumerate}
        \item \{1=23760, 2=596970, 3=3572100, 4=5354370\}\newline [[0, 1, 4, 55, 73, 121], [0, 2, 9, 95, 111, 128], [0, 3, 12, 29, 52, $\infty$], [0, 6, 54, 70, 116, 119],\newline   [0, 8, 15, 76, 89, 133], [0, 10, 71, 74, 107, 127]]
        \item \{1=25920, 2=635040, 3=3528360, 4=5357880\}\newline [[0, 1, 2, 30, 71, 106], [0, 3, 12, 29, 52, $\infty$], [0, 4, 31, 33, 46, 104], [0, 6, 10, 45, 84, 120],\newline  [0, 7, 24, 94, 111, 122], [0, 11, 40, 44, 48, 113]]
        \item \{1=27000, 2=627750, 3=3538620, 4=5353830\}\newline [[0, 1, 10, 50, 78, 124], [0, 2, 8, 9, 18, 42], [0, 3, 12, 29, 52, $\infty$], [0, 4, 88, 94, 114, 120],\newline  [0, 5, 91, 93, 106, 129], [0, 7, 39, 60, 68, 116]]
        \item \{1=29160, 2=671490, 3=3604500, 4=5242050\}\newline [[0, 1, 4, 24, 109, 129], [0, 3, 12, 29, 52, $\infty$], [0, 6, 77, 81, 87, 102], [0, 8, 21, 59, 103, 126],\newline [0, 11, 31, 39, 57, 66], [0, 25, 28, 71, 94, 113]]
        \item \{1=30240, 2=647190, 3=3546180, 4=5323590\}\newline [[0, 1, 2, 19, 71, 73], [0, 3, 12, 29, 52, $\infty$], [0, 4, 31, 67, 69, 124], [0, 5, 20, 47, 106, 108],\newline  [0, 6, 17, 43, 85, 134], [0, 11, 48, 89, 93, 96]]
        \item \{1=30240, 2=669060, 3=3601800, 4=5246100\}\newline [[0, 1, 4, 15, 96, 111], [0, 3, 12, 29, 52, $\infty$], [0, 7, 93, 99, 103, 134], [0, 8, 51, 66, 69, 132],\newline [0, 11, 46, 53, 92, 112], [0, 25, 28, 71, 94, 113]]
        \item \{1=32400, 2=631800, 3=3593160, 4=5289840\}\newline [[0, 1, 2, 44, 46, 134], [0, 3, 12, 29, 52, $\infty$], [0, 4, 24, 84, 102, 103], [0, 5, 7, 27, 95, 132],\newline  [0, 6, 17, 38, 80, 96], [0, 11, 48, 111, 114, 116]]
        \item \{1=32400, 2=633420, 3=3565080, 4=5316300\}\newline [[0, 1, 4, 27, 77, 122], [0, 2, 9, 19, 45, 64], [0, 3, 12, 29, 52, $\infty$], [0, 6, 49, 83, 99, 114],\newline  [0, 15, 16, 21, 25, 109], [0, 20, 70, 73, 94, 130]]
        \item \{1=34560, 2=647190, 3=3603420, 4=5262030\}\newline [[0, 1, 4, 22, 128, 130], [0, 3, 12, 29, 52, $\infty$], [0, 6, 63, 70, 80, 83], [0, 10, 21, 25, 53, 55],\newline  [0, 11, 48, 88, 108, 122], [0, 15, 31, 44, 50, 124]]
        \item \{1=34560, 2=663390, 3=3580740, 4=5268510\}\newline [[0, 1, 4, 16, 94, 114], [0, 3, 12, 29, 52, $\infty$], [0, 6, 65, 76, 95, 96], [0, 8, 24, 73, 85, 112],\newline  [0, 10, 63, 82, 123, 126], [0, 11, 48, 61, 84, 104]]
        \item \{1=35640, 2=656100, 3=3549960, 4=5305500\}\newline [[0, 1, 4, 10, 31, 106], [0, 3, 12, 29, 52, $\infty$], [0, 6, 22, 81, 93, 125], [0, 8, 82, 86, 89, 108],\newline  [0, 24, 30, 65, 66, 103], [0, 25, 28, 71, 94, 113]]
        \item \{1=36720, 2=648000, 3=3592080, 4=5270400\}\newline [[0, 1, 10, 56, 102, 117], [0, 3, 12, 29, 52, $\infty$], [0, 4, 63, 75, 81, 131], [0, 5, 18, 54, 93, 129],\newline   [0, 6, 26, 36, 59, 101], [0, 25, 28, 43, 66, 89]]
        \item \{1=36720, 2=660960, 3=3602880, 4=5246640\}\newline [[0, 1, 2, 19, 109, 118], [0, 3, 12, 29, 52, $\infty$], [0, 4, 24, 79, 93, 133], [0, 6, 46, 76, 103, 128],\newline [0, 8, 31, 60, 67, 132], [0, 11, 18, 20, 23, 48]]
        \item \{1=37800, 2=667440, 3=3589920, 4=5252040\}\newline [[0, 1, 4, 38, 45, 103], [0, 3, 12, 29, 52, $\infty$], [0, 6, 17, 80, 96, 107], [0, 7, 57, 65, 70, 134],\newline  [0, 10, 39, 93, 100, 111], [0, 11, 48, 62, 85, 105]]
        \item \{1=38880, 2=683640, 3=3591000, 4=5233680\}\newline [[0, 1, 2, 15, 38, 60], [0, 3, 12, 29, 52, $\infty$], [0, 4, 54, 86, 89, 117], [0, 5, 49, 50, 63, 102],\newline  [0, 11, 48, 111, 114, 116], [0, 19, 44, 55, 110, 129]]
        \item \{0=1350, 1=22680, 2=652860, 3=3636360, 4=5233950\}\newline [[0, 1, 2, 19, 68, 77], [0, 3, 12, 29, 52, $\infty$], [0, 5, 42, 81, 110, 129], [0, 6, 31, 71, 72, 131],\newline  [0, 11, 40, 44, 48, 113], [0, 13, 39, 46, 57, 95]]
        \item {\bf \{0=1350, 1=35640, 2=606690, 3=3615300, 4=5288220\}\newline [[0, 1, 4, 15, 51, 79], [0, 3, 12, 29, 52, $\infty$], [0, 8, 38, 77, 118, 124],\newline [0, 11, 48, 62, 85, 105], [0, 13, 66, 67, 75, 86], [0, 18, 41, 104, 106, 109]]}
        \item \{0=1350, 1=33480, 2=622890, 3=3562380, 4=5327100\}\newline [[0, 1, 10, 62, 66, 100], [0, 2, 8, 9, 18, 42], [0, 3, 12, 29, 52, $\infty$], [0, 4, 26, 99, 124, 130],\newline  [0, 5, 20, 31, 110, 134], [0, 7, 44, 50, 121, 132]]
        \item \{0=1350, 1=32400, 2=634230, 3=3591540, 4=5287680\}\newline [[0, 1, 2, 27, 51, 115], [0, 3, 12, 29, 52, $\infty$], [0, 5, 20, 22, 80, 133], [0, 6, 38, 53, 112, 134],\newline   [0, 7, 36, 42, 72, 120], [0, 11, 40, 44, 48, 113]]
        \item \{0=4050, 1=41040, 2=651240, 3=3621240, 4=5229630\}\newline [[0, 1, 2, 31, 67, 77], [0, 3, 12, 29, 52, $\infty$], [0, 4, 64, 70, 113, 121], [0, 5, 20, 73, 76, 134],\newline  [0, 6, 22, 30, 33, 120], [0, 11, 48, 111, 114, 116]]
    \end{enumerate}
    \end{example}

    \begin{example} There are at least four $1$-rotational designs for {\tt SmallGroup(150,7) = C2\,x\,((C5\,x\,C5)\,:\,C3)}
    \begin{enumerate}
        \item \{1=25200, 2=732600, 3=4572000, 4=7807200\} \newline [[0, 1, 2, 5, 8, 14], [0, 3, 9, 22, 114, 131], [0, 4, 53, 100, 134, 138], [0, 6, 10, 44, 50, 129], \newline [0, 11, 37, 52, 87, 98], [0, 20, 41, 65, 116, 119], [0, 31, 75, 94, 128, 130], \newline [0, 32, 86, 92, 102, 133], [0, 39, 70, 78, 105, 109], [0, 71, 72, 73, 74, $\infty$]]
        \item \{1=33600, 2=720000, 3=4539600, 4=7843800\} \newline [[0, 1, 2, 5, 8, 14], [0, 3, 9, 22, 114, 131], [0, 4, 17, 61, 121, 147], [0, 6, 23, 93, 133, 140], \newline [0, 10, 74, 78, 97, 136], [0, 18, 44, 90, 108, 137], [0, 27, 46, 99, 119, $\infty$], \newline [0, 29, 32, 59, 142, 144], [0, 30, 41, 81, 109, 115], [0, 35, 68, 100, 132, 139]]
        \item \{1=34800, 2=727200, 3=4588800, 4=7786200\} \newline [[0, 1, 2, 5, 8, 14], [0, 3, 9, 22, 114, 131], [0, 4, 75, 76, 98, 128], [0, 7, 35, 118, 132, 147], \newline [0, 12, 28, 33, 50, 141], [0, 19, 23, 64, 82, 124], [0, 20, 26, 69, 89, 116], \newline [0, 27, 46, 99, 119, $\infty$], [0, 30, 66, 108, 111, 113], [0, 32, 86, 92, 102, 133]]
        \item \{0=1500, 1=45600, 2=776700, 3=4596600, 4=7716600\} \newline [[0, 1, 2, 5, 8, 14], [0, 3, 9, 22, 114, 131], [0, 4, 80, 84, 97, 103], [0, 11, 28, 32, 81, 110], \newline [0, 18, 75, 85, 92, 116], [0, 19, 53, 89, 117, 137], [0, 23, 62, 64, 73, 134], \newline [0, 27, 46, 99, 119, $\infty$], [0, 35, 91, 95, 146, 147], [0, 36, 68, 101, 132, 140]]
    \end{enumerate}
    \end{example}

    \begin{example} There is at least one $1$-rotational design for {\tt SmallGroup(155,1) = C31\,:\,C5}
    \begin{enumerate}
        \item[] \{1=31000, 2=789570, 3=4994100, 4=8693330\}\newline [[0, 1, 4, 11, 13, 19], [0, 3, 22, 46, 127, 145], [0, 12, 30, 52, 104, 115], [0, 15, 21, 37, 73, $\infty$],\newline [0, 16, 32, 111, 112, 128], [0, 18, 72, 83, 94, 107]]
    \end{enumerate}
    \end{example}

    \begin{example} There are at least sixteen $1$-rotational designs for the cyclic  group $\mathbb Z_{155}$.
    \begin{enumerate}
        \item \{1=22320, 2=771900, 3=4973640, 4=8740140\}\newline [[0, 1, 3, 24, 40, 73], [0, 4, 13, 78, 100, 108], [0, 5, 61, 102, 128, 145],  [0, 6, 50, 107, 119, 126], [0, 11, 25, 45, 63, 91], [0, 31, 62, 93, 124, $\infty$]]
        \item \{1=24800, 2=720750, 3=4997820, 4=8764630\}\newline [[0, 1, 3, 24, 40, 73], [0, 4, 13, 78, 100, 108], [0, 5, 61, 102, 128, 145],  [0, 6, 50, 107, 119, 126], [0, 11, 75, 103, 121, 141], [0, 31, 62, 93, 124, $\infty$]]
        \item \{1=26040, 2=716100, 3=4973640, 4=8792220\}\newline [[0, 1, 3, 24, 40, 73], [0, 4, 13, 78, 100, 108], [0, 5, 15, 32, 58, 99], [0, 6, 50, 107, 119, 126],\newline   [0, 11, 25, 45, 63, 91], [0, 31, 62, 93, 124, $\infty$]]
        \item \{1=28520, 2=762600, 3=4910400, 4=8806480\}\newline [[0, 1, 3, 24, 40, 73], [0, 4, 51, 59, 81, 146], [0, 5, 15, 32, 58, 99], [0, 6, 50, 107, 119, 126],\newline   [0, 11, 25, 45, 63, 91], [0, 31, 62, 93, 124, $\infty$]]
        \item \{1=29760, 2=724470, 3=4990380, 4=8763390\}\newline [[0, 1, 3, 24, 40, 73], [0, 4, 13, 78, 100, 108], [0, 5, 15, 32, 58, 99], [0, 6, 50, 107, 119, 126],\newline   [0, 11, 75, 103, 121, 141], [0, 31, 62, 93, 124, $\infty$]]
        \item \{1=31000, 2=772830, 3=4997820, 4=8706350\}\newline [[0, 1, 3, 24, 40, 73], [0, 4, 51, 59, 81, 146], [0, 5, 15, 32, 58, 99], [0, 6, 35, 42, 54, 111],\newline  [0, 11, 25, 45, 63, 91], [0, 31, 62, 93, 124, $\infty$]]
        \item \{1=31000, 2=823980, 3=5022000, 4=8631020\}\newline [[0, 1, 3, 24, 40, 73], [0, 4, 51, 59, 81, 146], [0, 5, 15, 32, 58, 99], [0, 6, 35, 42, 54, 111],\newline [0, 11, 75, 103, 121, 141], [0, 31, 62, 93, 124, $\infty$]]
        \item \{1=33480, 2=829560, 3=4992240, 4=8652720\}\newline [[0, 1, 3, 24, 40, 73], [0, 4, 51, 59, 81, 146], [0, 5, 61, 102, 128, 145], [0, 6, 50, 107, 119, 126], [0, 11, 75, 103, 121, 141], [0, 31, 62, 93, 124, $\infty$]]
        \item \{1=34720, 2=839790, 3=4960620, 4=8672870\}\newline [[0, 1, 3, 24, 40, 73], [0, 4, 51, 59, 81, 146], [0, 5, 61, 102, 128, 145], [0, 6, 50, 107, 119, 126], [0, 11, 25, 45, 63, 91], [0, 31, 62, 93, 124, $\infty$]]
        \item \{1=35960, 2=817470, 3=5068500, 4=8586070\}\newline [[0, 1, 3, 24, 40, 73], [0, 4, 13, 78, 100, 108], [0, 5, 15, 32, 58, 99], [0, 6, 35, 42, 54, 111],\newline  [0, 11, 75, 103, 121, 141], [0, 31, 62, 93, 124, $\infty$]]
        \item \{1=42160, 2=746790, 3=4960620, 4=8758430\}\newline [[0, 1, 3, 24, 40, 73], [0, 4, 13, 78, 100, 108], [0, 5, 61, 102, 128, 145], [0, 6, 35, 42, 54, 111], [0, 11, 75, 103, 121, 141], [0, 31, 62, 93, 124, $\infty$]]
        \item \{1=42160, 2=795150, 3=4997820, 4=8672870\}\newline [[0, 1, 3, 24, 40, 73], [0, 4, 51, 59, 81, 146], [0, 5, 61, 102, 128, 145], [0, 6, 35, 42, 54, 111],\newline   [0, 11, 75, 103, 121, 141], [0, 31, 62, 93, 124, $\infty$]]
        \item \{1=44640, 2=763530, 3=4934580, 4=8765250\}\newline [[0, 1, 3, 24, 40, 73], [0, 4, 51, 59, 81, 146], [0, 5, 15, 32, 58, 99], [0, 6, 50, 107, 119, 126],\newline   [0, 11, 75, 103, 121, 141], [0, 31, 62, 93, 124, $\infty$]]
        \item \{0=1550, 1=32240, 2=821190, 3=5021380, 4=8631640\}\newline [[0, 1, 3, 24, 40, 73], [0, 4, 13, 78, 100, 108], [0, 5, 61, 102, 128, 145], [0, 6, 35, 42, 54, 111], [0, 11, 25, 45, 63, 91], [0, 31, 62, 93, 124, $\infty$]]
        \item \{0=1550, 1=42160, 2=812820, 3=5026960, 4=8624510\}\newline [[0, 1, 3, 24, 40, 73], [0, 4, 51, 59, 81, 146], [0, 5, 61, 102, 128, 145], [0, 6, 35, 42, 54, 111],\newline   [0, 11, 25, 45, 63, 91], [0, 31, 62, 93, 124, $\infty$]]
        \item \{0=1550, 1=43400, 2=788640, 3=5038120, 4=8636290\}\newline [[0, 1, 3, 24, 40, 73], [0, 4, 13, 78, 100, 108], [0, 5, 15, 32, 58, 99], [0, 6, 35, 42, 54, 111],\newline  [0, 11, 25, 45, 63, 91], [0, 31, 62, 93, 124, $\infty$]]
    \end{enumerate}
    \end{example}

    \section{New $S(2,6,66)$ with automorphism group of order 20}

    In this section we will describe new Steiner system $S(2,6,66)$ generated by group $\mathbb Z_2 \times \mathbb Z_{10}$ with five orbits of size $20, 20, 20, 5, 1$. Elements of first three orbits will be denoted by numbers $0',\dots,19'$, $0'',\dots,19''$ and $0''',\dots,19'''$ with one, two and three primes respectively. Elements of the $5$-element orbit is denoted by numbers $0,\dots,5$, and element from the one-element orbit is denoted by $\infty$. Group action on the $20$-element orbits are defined by the formula
    $$x \oplus y =
    \begin{cases}
    (x + y) \mod 10, &\mbox{if $(x < 10 \land y < 10) \lor (x \ge 10 \land y \ge 10)$};\\
    10 + ((x + y) \mod 10), &\mbox{otherwise}.
    \end{cases}
    $$
    The group action on the 5-element orbit is defined as $x \oplus y = (x + y) \mod 5$.

    \begin{example} Base blocks of design are:\newline
    [[$0'$, $2'$, $13'$, $17''$, $11'''$, $12'''$], [$0'$, $5'$, $10'$, $15'$, $0$, $\infty$],  [$0'$, $6'$, $12'$, $12''$, $19''$, $0'''$],\newline [$0'$, $1''$, $2''$, $18''$, $3'''$, $15'''$], [$0'$, $5''$, $14''$, $1'''$, $14'''$, $17'''$],  [$2'$, $3'$, $1''$, $13''$, $9'''$, $0$],\newline [$6'$, $9'$, $2''$, $1'''$, $12'''$, $0$], [$0''$, $2''$, $4''$, $6''$, $8''$, $\infty$], [$0''$, $5''$, $10''$, $15''$, $0$, $1$],\newline  [$0'''$, $2'''$, $4'''$, $6'''$, $8'''$, $\infty$], [$0'''$, $5'''$, $10'''$, $15'''$, $0$, $2$]]
    \end{example}

    \section{Designs generated by a cyclic group and two orbits}

    In this section we present Steiner systems $S(2,6,2n)$ for $n\in\{48,53\}$ which are generated by the action of the cyclic group $\mathbb Z_n$ on two orbits of length $n$. Elements of the first orbit are identified with elements $0,\dots,n-1$ of the cyclic groups $\mathbb Z_n$ and the elements of the second orbit are labeled by the numbers $0',\dots,(n-1)'$ endowed with primes.

    \begin{example} There exists only one design for $v=96$ generated by the cyclic group $\mathbb Z_{48}$ acting on two orbits of size $48$. It is Mills' design \cite{Mills1}.
    \begin{enumerate}
        \item \{1=37632, 2=492192, 3=1475136, 4=1278240\}\newline [[$0, 8, 16, 24, 32, 40$], [$0', 8', 16', 24', 32', 40'$], [$0, 1, 3, 13, 28, 0'$], [$0, 4, 11, 17', 36', 38'$],\newline [$0, 5, 19, 1', 24', 42'$], [$0, 9, 26, 4', 7', 40'$], [$0, 6, 8', 9', 18', 22'$], [$0, 18, 11', 28', 33', 39'$]]
    \end{enumerate}
    \end{example}

    \begin{example} There exist 66 designs for $v=106$ generated by the group $\mathbb Z_{53}$ acting on two orbits of size $53$. The design (46) was foubd by Mills  \cite{Mills2}.
    \begin{enumerate}
        \item \{1=24168, 2=491628, 3=1965240, 4=1970964\}\newline [[$0, 1, 4, 11, 0', 19'$], [$0, 2, 24, 37, 7', 12'$], [$0, 5, 30, 39, 2', 3'$], [$0, 6, 26, 38, 30', 39'$],\newline  [$0, 8, 14', 29', 43', 46'$], [$0, 17, 11', 37', 44', 48'$], [$0, 9', 22', 32', 34', 40'$]]
        \item \{1=27136, 2=504984, 3=1948704, 4=1971176\}\newline [[$0, 1, 3, 39, 0', 6'$], [$0, 4, 28, 41, 46, 22'$], [$0, 6, 32, 30', 43', 44'$], [$0, 8, 31, 39', 41', 48'$],\newline  [$0, 9, 19, 2', 28', 32'$], [$0, 20, 16', 27', 35', 45'$], [$0, 1', 4', 21', 26', 42'$]]
        \item \{1=29256, 2=490356, 3=1974144, 4=1958244\}\newline [[$0, 1, 3, 12, 40, 0'$], [$0, 4, 30, 36, 37', 48'$], [$0, 5, 20, 26', 28', 34'$], [$0, 7, 31, 9', 10', 47'$],\newline  [$0, 8, 18, 4', 38', 43'$], [$0, 19, 11', 24', 36', 46'$], [$0, 15', 19', 22', 42', 51'$]]
        \item \{1=31800, 2=505302, 3=1977324, 4=1937574\}\newline [[$0, 1, 3, 10, 31, 0'$], [$0, 4, 12, 18, 21', 37'$], [$0, 5, 24, 6', 13', 28'$], [$0, 11, 38, 5', 45', 49'$],\newline  [$0, 13, 15', 29', 39', 40'$], [$0, 16, 36, 14', 46', 48'$], [$0, 18', 24', 36', 41', 44'$]]
        \item \{1=32224, 2=524382, 3=1944252, 4=1951142\}\newline [[$0, 1, 3, 15, 35, 0'$], [$0, 4, 29, 27', 36', 43'$], [$0, 5, 31, 9', 26', 47'$], [$0, 6, 16, 46, 8', 12'$],\newline  [$0, 8, 44, 1', 25', 28'$], [$0, 11, 22', 35', 40', 41'$], [$0, 3', 5', 13', 33', 44'$]]
        \item \{1=32648, 2=504666, 3=1959516, 4=1955170\}\newline [[$0, 1, 3, 13, 49, 0'$], [$0, 6, 28, 11', 37', 48'$], [$0, 8, 27, 15', 18', 33'$], [$0, 9, 24, 42, 1', 32'$],\newline  [$0, 14, 30, 16', 28', 49'$], [$0, 21, 24', 34', 38', 47'$], [$0, 21', 22', 27', 29', 46'$]]
        \item \{1=32648, 2=512298, 3=1979868, 4=1927186\}\newline [[$0, 1, 3, 10, 0', 12'$], [$0, 4, 36, 14', 25', 32'$], [$0, 5, 30, 6', 35', 45'$], [$0, 6, 18, 26, 37, 44'$],\newline  [$0, 13, 17', 33', 36', 37'$], [$0, 14, 29, 22', 27', 48'$], [$0, 3', 16', 39', 41', 47'$]]
        \item \{1=33072, 2=506574, 3=1965876, 4=1946478\}\newline [[$0, 1, 3, 15, 31, 0'$], [$0, 4, 36, 44, 31', 45'$], [$0, 5, 20', 33', 39', 49'$], [$0, 6, 33, 10', 17', 35'$],\newline  [$0, 7, 18, 23', 25', 26'$], [$0, 10, 34, 13', 24', 46'$], [$0, 6', 21', 42', 47', 51'$]]
        \item \{1=34344, 2=528516, 3=1938528, 4=1950612\}\newline [[$0, 1, 3, 8, 23, 0'$], [$0, 4, 16, 40, 19', 44'$], [$0, 6, 32, 18', 27', 37'$], [$0, 9, 19, 17', 33', 35'$],\newline  [$0, 11, 25, 1', 34', 47'$], [$0, 18, 6', 7', 11', 38'$], [$0, 2', 10', 13', 25', 49'$]]
        \item \{1=34768, 2=506574, 3=1958244, 4=1952414\}\newline [[$0, 1, 3, 15, 45, 0'$], [$0, 4, 26, 37', 43', 44'$], [$0, 5, 25, 29', 31', 46'$], [$0, 6, 24, 40, 7', 34'$],\newline  [$0, 7, 43, 3', 12', 32'$], [$0, 21, 16', 19', 30', 35'$], [$0, 2', 15', 23', 27', 45'$]]
        \item \{1=35616, 2=513570, 3=1979868, 4=1922946\}\newline [[$0, 1, 3, 30, 38, 0'$], [$0, 4, 13, 16', 22', 46'$], [$0, 5, 11, 25, 19', 35'$], [$0, 7, 41, 28', 43', 45'$],\newline  [$0, 10, 32, 6', 11', 39'$], [$0, 17, 5', 37', 48', 51'$], [$0, 13', 17', 25', 26', 44'$]]
        \item \{1=35616, 2=521838, 3=1970964, 4=1923582\}\newline [[$0, 1, 3, 15, 0', 9'$], [$0, 4, 13, 20, 30', 41'$], [$0, 5, 24, 47, 25', 29'$], [$0, 8, 26, 36, 15', 22'$],\newline  [$0, 21, 2', 4', 16', 19'$], [$0, 22, 12', 13', 40', 45'$], [$0, 3', 11', 27', 33', 46'$]]
        \item \{1=36040, 2=492900, 3=1967784, 4=1955276\}\newline [[$0, 1, 4, 15, 44, 0'$], [$0, 2, 18, 4', 26', 36'$], [$0, 5, 26, 46, 33', 51'$], [$0, 6, 34, 3', 16', 23'$],\newline  [$0, 8, 31, 19', 21', 45'$], [$0, 17, 32', 44', 47', 48'$], [$0, 1', 6', 12', 20', 29'$]]
        \item \{1=36464, 2=489084, 3=1989408, 4=1937044\}\newline [[$0, 1, 3, 9, 0', 7'$], [$0, 4, 16, 26, 40, 49'$], [$0, 5, 18', 29', 32', 42'$], [$0, 7, 32, 10', 22', 48'$],\newline  [$0, 11, 34, 19', 28', 36'$], [$0, 15, 33, 1', 20', 26'$], [$0, 12', 14', 30', 34', 35'$]]
        \item \{1=36464, 2=516432, 3=1971600, 4=1927504\}\newline [[$0, 1, 3, 39, 0', 47'$], [$0, 4, 32, 45, 16', 21'$], [$0, 5, 11, 35, 15', 41'$], [$0, 7, 16, 26, 35', 39'$],\newline  [$0, 20, 7', 22', 25', 38'$], [$0, 22, 3', 20', 48', 49'$], [$0, 1', 11', 31', 43', 45'$]]
        \item \{1=37312, 2=501804, 3=1942344, 4=1970540\}\newline [[$0, 1, 3, 13, 21, 0'$], [$0, 4, 11', 28', 42', 49'$], [$0, 5, 19, 28, 14', 41'$], [$0, 6, 22, 25', 26', 37'$],\newline  [$0, 7, 36, 12', 17', 30'$], [$0, 11, 38, 1', 44', 46'$], [$0, 2', 18', 21', 27', 51'$]]
        \item \{1=37736, 2=502122, 3=1963332, 4=1948810\}\newline [[$0, 1, 3, 21, 46, 0'$], [$0, 4, 19, 24', 33', 46'$], [$0, 5, 31, 47, 17', 41'$], [$0, 9, 39, 30', 34', 35'$],\newline  [$0, 12, 29, 13', 16', 31'$], [$0, 13, 3', 11', 22', 28'$], [$0, 6', 8', 18', 38', 45'$]]
        \item \{1=37736, 2=530424, 3=1941072, 4=1942768\}\newline [[$0, 1, 3, 23, 41, 0'$], [$0, 4, 28, 47, 2', 42'$], [$0, 5, 6', 15', 16', 39'$], [$0, 7, 21, 28', 40', 44'$],\newline  [$0, 8, 44, 4', 26', 32'$], [$0, 11, 27, 3', 5', 20'$], [$0, 17', 22', 25', 36', 43'$]]
        \item \{1=38160, 2=493854, 3=1970964, 4=1949022\}\newline [[$0, 1, 3, 14, 0', 35'$], [$0, 4, 27, 36, 5', 19'$], [$0, 5, 24, 46, 30', 33'$], [$0, 6, 4', 20', 48', 49'$],\newline  [$0, 8, 28, 43, 16', 46'$], [$0, 16, 10', 23', 27', 29'$], [$0, 2', 12', 17', 24', 44'$]]
        \item \{1=38160, 2=512934, 3=1950612, 4=1950294\}\newline [[$0, 1, 3, 23, 0', 6'$], [$0, 4, 10, 14', 21', 43'$], [$0, 5, 13, 32, 41, 20'$], [$0, 7, 42, 2', 23', 35'$],\newline  [$0, 14, 29, 1', 38', 51'$], [$0, 16, 8', 12', 42', 47'$], [$0, 18', 19', 27', 29', 44'$]]
        \item \{1=39008, 2=491628, 3=1938528, 4=1982836\}\newline [[$0, 1, 3, 7, 22, 0'$], [$0, 5, 30, 42, 22', 48'$], [$0, 8, 43, 4', 13', 34'$], [$0, 9, 29, 38', 39', 41'$],\newline  [$0, 13, 27, 11', 21', 28'$], [$0, 17, 20', 24', 36', 42'$], [$0, 2', 16', 27', 35', 40'$]]
        \item \{1=39432, 2=529152, 3=1930896, 4=1952520\}\newline [[$0, 1, 3, 30, 48, 0'$], [$0, 4, 16, 24', 25', 45'$], [$0, 7, 21, 43, 22', 34'$], [$0, 9, 20, 3', 16', 39'$],\newline  [$0, 13, 38, 2', 11', 48'$], [$0, 19, 6', 31', 33', 37'$], [$0, 4', 28', 38', 43', 46'$]]
        \item \{1=39856, 2=501804, 3=1972872, 4=1937468\}\newline [[$0, 1, 3, 15, 46, 0'$], [$0, 4, 36, 16', 24', 27'$], [$0, 5, 16, 40, 6', 26'$], [$0, 6, 33, 2', 31', 46'$],\newline  [$0, 9, 28, 3', 4', 17'$], [$0, 23, 5', 15', 32', 37'$], [$0, 11', 18', 30', 34', 36'$]]
        \item \{1=40280, 2=542190, 3=1942980, 4=1926550\}\newline [[$0, 1, 3, 30, 44, 0'$], [$0, 4, 25, 7', 10', 19'$], [$0, 5, 22, 38, 33', 39'$], [$0, 6, 19, 24', 46', 51'$],\newline  [$0, 7, 18, 20', 43', 44'$], [$0, 8, 12', 16', 30', 49'$], [$0, 14', 21', 29', 31', 42'$]]
        \item \{1=41976, 2=523746, 3=1959516, 4=1926762\}\newline [[$0, 1, 3, 9, 40, 0'$], [$0, 4, 32, 9', 14', 28'$], [$0, 5, 17, 23', 25', 51'$], [$0, 7, 18, 33, 22', 45'$],\newline  [$0, 10, 34, 21', 29', 41'$], [$0, 23, 2', 17', 26', 39'$], [$0, 1', 33', 36', 37', 43'$]]
        \item \{1=41976, 2=529788, 3=1966512, 4=1913724\}\newline [[$0, 1, 3, 8, 36, 0'$], [$0, 4, 43, 8', 20', 39'$], [$0, 6, 27, 15', 25', 28'$], [$0, 9, 22, 38, 11', 43'$],\newline  [$0, 11, 30, 14', 23', 40'$], [$0, 12, 6', 7', 36', 44'$], [$0, 13', 27', 31', 33', 38'$]]
        \item \{1=41976, 2=536784, 3=1938528, 4=1934712\}\newline [[$0, 1, 3, 16, 0', 18'$], [$0, 4, 21, 26, 46, 29'$], [$0, 6, 45, 16', 19', 40'$], [$0, 9, 5', 31', 35', 41'$],\newline  [$0, 10, 34, 1', 38', 43'$], [$0, 12, 30, 23', 42', 51'$], [$0, 6', 7', 14', 45', 47'$]]
        \item \{1=42400, 2=509118, 3=1964604, 4=1935878\}\newline [[$0, 1, 3, 8, 21, 0'$], [$0, 4, 47, 6', 22', 33'$], [$0, 9, 26, 37, 10', 24'$], [$0, 12, 31, 16', 21', 25'$],\newline  [$0, 14, 29, 3', 48', 49'$], [$0, 23, 11', 14', 31', 46'$], [$0, 5', 7', 17', 30', 36'$]]
        \item \{1=43672, 2=517386, 3=1954428, 4=1936514\}\newline [[$0, 1, 3, 8, 29, 0'$], [$0, 4, 20, 38, 14', 39'$], [$0, 6, 17, 4', 33', 42'$], [$0, 9, 23, 11', 31', 32'$],\newline  [$0, 10, 22, 6', 13', 48'$], [$0, 13, 18', 28', 30', 34'$], [$0, 7', 12', 20', 43', 46'$]]
        \item \{1=44096, 2=499896, 3=1952520, 4=1955488\}\newline [[$0, 1, 3, 10, 41, 0'$], [$0, 4, 20, 39, 18', 30'$], [$0, 5, 26, 1', 7', 29'$], [$0, 6, 29, 15', 23', 33'$],\newline  [$0, 8, 25, 13', 45', 46'$], [$0, 11, 6', 19', 22', 36'$], [$0, 16', 31', 35', 40', 42'$]]
        \item \{1=44096, 2=529152, 3=1962696, 4=1916056\}\newline [[$0, 1, 3, 23, 41, 0'$], [$0, 4, 21, 29, 25', 31'$], [$0, 5, 11', 13', 22', 40'$], [$0, 6, 9', 29', 45', 48'$],\newline  [$0, 7, 16, 26, 14', 44'$], [$0, 11, 1', 16', 26', 47'$], [$0, 14, 33', 34', 38', 46'$]]
        \item \{1=47064, 2=543780, 3=1942344, 4=1918812\}\newline [[$0, 1, 3, 8, 19, 0'$], [$0, 4, 21, 44, 5', 23'$], [$0, 6, 33, 4', 28', 42'$], [$0, 10, 41, 13', 17', 26'$],\newline  [$0, 14, 29, 41', 44', 49'$], [$0, 25, 11', 31', 33', 43'$], [$0, 10', 21', 40', 46', 47'$]]
        \item \{1=48336, 2=530424, 3=1957608, 4=1915632\}\newline [[$0, 1, 3, 22, 46, 0'$], [$0, 4, 16, 39, 22', 48'$], [$0, 5, 11, 21', 38', 40'$], [$0, 9, 36, 3', 28', 51'$],\newline  [$0, 13, 28, 14', 24', 30'$], [$0, 20, 4', 8', 43', 46'$], [$0, 5', 12', 13', 25', 34'$]]

        \item \{0=530, 1=33496, 2=509436, 3=1950400, 4=1958138\}\newline [[$0, 1, 4, 17, 0', 2'$], [$0, 2, 21, 30, 44', 48'$], [$0, 5, 31, 46, 8', 9'$], [$0, 6, 14, 24, 34', 43'$],\newline  [$0, 11, 22', 32', 35', 50'$], [$0, 20, 7', 12', 26', 33'$], [$0, 5', 17', 25', 41', 47'$]]
        \item \{0=530, 1=33920, 2=505302, 3=1967572, 4=1944676\}\newline [[$0, 1, 3, 14, 36, 0'$], [$0, 4, 45, 7', 25', 36'$], [$0, 5, 32, 28', 43', 51'$], [$0, 6, 34, 44, 1', 18'$],\newline  [$0, 7, 37, 29', 31', 41'$], [$0, 24, 6', 13', 26', 40'$], [$0, 5', 8', 9', 14', 30'$]]
        \item \{0=530, 1=34344, 2=534240, 3=1956760, 4=1926126\}\newline [[$0, 1, 3, 13, 35, 0'$], [$0, 4, 48, 28', 42', 43'$], [$0, 6, 14, 29, 17', 36'$], [$0, 7, 27, 8', 19', 32'$],\newline  [$0, 11, 36, 4', 20', 27'$], [$0, 16, 26', 29', 31', 51'$], [$0, 2', 6', 14', 23', 49'$]]
        \item \{0=530, 1=35616, 2=499896, 3=1972024, 4=1943934\}\newline [[$0, 1, 3, 7, 0', 12'$], [$0, 5, 23, 32, 42, 49'$], [$0, 8, 22, 37', 40', 41'$], [$0, 12, 36, 13', 22', 28'$],\newline  [$0, 13, 28, 2', 34', 48'$], [$0, 20, 14', 24', 43', 51'$], [$0, 3', 8', 25', 36', 38'$]]
        \item \{0=530, 1=35616, 2=512616, 3=1960576, 4=1942662\}\newline [[$0, 1, 3, 9, 42, 0'$], [$0, 4, 25, 38, 12', 14'$], [$0, 5, 23, 18', 26', 51'$], [$0, 7, 31, 16', 32', 37'$],\newline  [$0, 10, 37, 4', 15', 33'$], [$0, 17, 7', 19', 34', 41'$], [$0, 22', 35', 36', 39', 45'$]]
        \item \{0=530, 1=36040, 2=516114, 3=1940860, 4=1958456\}\newline [[$0, 1, 3, 33, 39, 0'$], [$0, 4, 13, 9', 11', 46'$], [$0, 5, 27, 46, 15', 30'$], [$0, 8, 24, 2', 16', 43'$],\newline  [$0, 10, 28, 1', 23', 34'$], [$0, 11, 29', 32', 38', 39'$], [$0, 4', 12', 17', 36', 40'$]]
        \item \{0=530, 1=36888, 2=506574, 3=1954852, 4=1953156\}\newline [[$0, 1, 3, 11, 30, 0'$], [$0, 4, 20, 11', 29', 32'$], [$0, 5, 41, 27', 35', 46'$], [$0, 6, 31, 46, 14', 26'$],\newline  [$0, 9, 13', 15', 19', 43'$], [$0, 14, 35, 31', 38', 51'$], [$0, 1', 2', 18', 40', 45'$]]
        \item \{0=530, 1=36888, 2=529470, 3=1957396, 4=1927716\}\newline [[$0, 1, 3, 24, 38, 0'$], [$0, 4, 26, 7', 9', 12'$], [$0, 5, 12, 1', 28', 45'$], [$0, 6, 40, 19', 30', 31'$],\newline  [$0, 8, 33, 44, 2', 26'$], [$0, 10, 14', 20', 27', 48'$], [$0, 6', 21', 37', 41', 51'$]]
        \item \{0=530, 1=37312, 2=490038, 3=1981564, 4=1942556\}\newline [[$0, 1, 3, 14, 36, 0'$], [$0, 4, 45, 10', 11', 24'$], [$0, 5, 28, 30', 36', 40'$], [$0, 6, 38, 1', 29', 51'$],\newline  [$0, 7, 34, 44, 28', 49'$], [$0, 24, 4', 9', 27', 46'$], [$0, 14', 26', 34', 41', 43'$]]
        \item \{0=530, 1=38160, 2=554592, 3=1933864, 4=1924854\}\newline [[$0, 1, 3, 16, 0', 44'$], [$0, 4, 10, 21, 33, 46'$], [$0, 5, 10', 26', 29', 40'$], [$0, 7, 25, 34, 3', 11'$],\newline  [$0, 8, 9', 14', 15', 27'$], [$0, 14, 8', 12', 32', 34'$], [$0, 22, 2', 17', 38', 45'$]]
        \item \{0=530, 1=38584, 2=527244, 3=1977112, 4=1908530\}\newline [[$0, 1, 3, 13, 38, 0'$], [$0, 4, 21, 30, 9', 12'$], [$0, 5, 47, 21', 39', 43'$], [$0, 7, 29, 2', 13', 36'$],\newline  [$0, 8, 1', 18', 28', 30'$], [$0, 14, 33, 3', 31', 47'$], [$0, 4', 11', 19', 24', 25'$]]
        \item \{0=530, 1=39008, 2=514206, 3=1972660, 4=1925596\}\newline [[$0, 1, 3, 12, 17, 0'$], [$0, 4, 22, 11', 16', 23'$], [$0, 6, 21, 46, 30', 49'$], [$0, 8, 27, 13', 33', 48'$],\newline  [$0, 10, 30, 8', 32', 45'$], [$0, 24, 17', 34', 38', 44'$], [$0, 4', 18', 26', 27', 29'$]]
        \item {\bf \{0=530, 1=41128, 2=521202, 3=1957396, 4=1931744\}}\newline $\mathbf{[[0, 1, 3, 11, 38, 0'],\; [0, 4, 28', 30', 37', 47'],\; [0, 5, 19, 25, 36', 39'],\; [0, 7, 29, 8', 16', 48']}$,\newline  $\mathbf{[0, 9, 21, 12', 13', 27'],\; [0, 13, 30, 23', 35', 51'],\; [0, 2', 7', 25', 29', 49']]}$
        \item \{0=530, 1=41552, 2=520566, 3=1943404, 4=1945948\}\newline [[$0, 1, 3, 8, 35, 0'$], [$0, 4, 14, 27', 44', 46'$], [$0, 6, 14', 25', 28', 41'$], [$0, 9, 29, 5', 12', 33'$],\newline  [$0, 11, 28, 41, 26', 48'$], [$0, 15, 31, 9', 17', 21'$], [$0, 1', 10', 11', 16', 34'$]]
        \item \{0=530, 1=41976, 2=531696, 3=1933864, 4=1943934\}\newline [[$0, 1, 3, 49, 0', 24'$], [$0, 6, 21, 34, 12', 32'$], [$0, 8, 18, 30, 43', 45'$], [$0, 9, 20, 36, 29', 39'$],\newline  [$0, 14, 1', 8', 16', 48'$], [$0, 24, 7', 38', 41', 42'$], [$0, 5', 10', 22', 33', 49'$]]
        \item \{0=530, 1=41976, 2=532650, 3=1957396, 4=1919448\}\newline [[$0, 1, 3, 19, 24, 0'$], [$0, 4, 47, 20', 39', 51'$], [$0, 7, 15, 27, 17', 30'$], [$0, 9, 22, 28', 31', 33'$],\newline  [$0, 11, 36, 25', 48', 49'$], [$0, 14, 5', 21', 32', 41'$], [$0, 1', 8', 36', 40', 46'$]]
        \item \{0=530, 1=42400, 2=516114, 3=1963756, 4=1929200\}\newline [[$0, 1, 3, 13, 17, 0'$], [$0, 5, 24, 30, 8', 34'$], [$0, 7, 15, 33, 22', 26'$], [$0, 9, 27', 41', 47', 48'$],\newline  [$0, 11, 13', 16', 25', 35'$], [$0, 21, 1', 12', 17', 30'$], [$0, 22, 20', 28', 43', 45'$]]
        \item \{0=530, 1=42824, 2=518976, 3=1963120, 4=1926550\}\newline [[$0, 1, 3, 15, 48, 0'$], [$0, 4, 30, 7', 28', 43'$], [$0, 7, 16, 35, 17', 41'$], [$0, 10, 32, 12', 26', 46'$],\newline  [$0, 11, 24, 15', 19', 42'$], [$0, 17, 9', 37', 40', 49'$], [$0, 11', 21', 22', 27', 29'$]]
        \item \{0=530, 1=44520, 2=512934, 3=1971388, 4=1922628\}\newline [[$0, 1, 3, 7, 43, 0'$], [$0, 5, 26, 34, 32', 39'$], [$0, 9, 23, 18', 35', 47'$], [$0, 12, 37, 1', 15', 28'$],\newline  [$0, 15, 35, 11', 19', 22'$], [$0, 22, 2', 30', 36', 45'$], [$0, 20', 21', 25', 41', 43'$]]
        \item \{0=530, 1=47912, 2=542826, 3=1959940, 4=1900792\}\newline [[$0, 1, 3, 15, 46, 0'$], [$0, 4, 23, 48, 14', 29'$], [$0, 6, 32, 30', 33', 43'$], [$0, 11, 29, 16', 42', 47'$],\newline  [$0, 13, 15', 17', 21', 35'$], [$0, 16, 36, 9', 28', 39'$], [$0, 20', 32', 41', 48', 49'$]]
        \item \{0=1060, 1=30528, 2=503394, 3=1966724, 4=1950294\}\newline [[$0, 1, 3, 14, 22, 0'$], [$0, 4, 10, 8', 17', 29'$], [$0, 5, 30, 46, 20', 26'$], [$0, 9, 26, 3', 32', 37'$],\newline  [$0, 15, 33, 2', 16', 24'$], [$0, 24, 5', 12', 38', 42'$], [$0, 10', 35', 45', 46', 48'$]]
        \item \{0=1060, 1=31376, 2=492264, 3=1959728, 4=1967572\}\newline [[$0, 1, 3, 7, 36, 0'$], [$0, 5, 45, 16', 40', 47'$], [$0, 9, 32, 43, 23', 41'$], [$0, 12, 39, 20', 25', 31'$],\newline  [$0, 15, 31, 27', 36', 37'$], [$0, 25, 10', 26', 29', 43'$], [$0, 3', 7', 15', 28', 30'$]]
        \item \{0=1060, 1=31800, 2=507846, 3=1967996, 4=1943298\}\newline [[$0, 1, 3, 15, 0', 46'$], [$0, 4, 21, 28, 48, 30'$], [$0, 6, 43, 12', 13', 48'$], [$0, 8, 42, 16', 25', 28'$],\newline  [$0, 13, 31, 11', 32', 34'$], [$0, 23, 10', 37', 41', 47'$], [$0, 4', 15', 29', 44', 49'$]]
        \item \{0=1060, 1=32648, 2=519612, 3=1974992, 4=1923688\}\newline [[$0, 1, 3, 7, 38, 0'$], [$0, 5, 34, 1', 29', 47'$], [$0, 8, 25, 10', 16', 36'$], [$0, 9, 30, 42, 12', 14'$],\newline  [$0, 10, 27', 32', 43', 51'$], [$0, 13, 27, 19', 31', 34'$], [$0, 9', 26', 30', 39', 40'$]]
        \item \{0=1060, 1=34344, 2=495762, 3=1959092, 4=1961742\}\newline [[$0, 1, 3, 7, 34, 0'$], [$0, 5, 17, 20', 31', 40'$], [$0, 8, 43, 6', 32', 45'$], [$0, 9, 32, 17', 22', 39'$],\newline  [$0, 11, 25, 40, 5', 21'$], [$0, 16, 4', 25', 27', 28'$], [$0, 1', 29', 36', 44', 48'$]]
        \item \{0=1060, 1=39008, 2=520248, 3=1962272, 4=1929412\}\newline [[$0, 1, 3, 8, 23, 0'$], [$0, 4, 29, 39, 5', 11'$], [$0, 6, 40, 15', 38', 43'$], [$0, 9, 41, 22', 36', 49'$],\newline  [$0, 11, 37, 2', 23', 31'$], [$0, 17, 6', 10', 21', 41'$], [$0, 14', 16', 17', 26', 33'$]]
        \item \{0=1060, 1=40280, 2=524382, 3=1952732, 4=1933546\}\newline [[$0, 1, 4, 22, 24, 0'$], [$0, 5, 14, 41, 3', 23'$], [$0, 6, 43, 16', 27', 34'$], [$0, 7, 2', 4', 40', 43'$],\newline  [$0, 8, 42, 1', 14', 30'$], [$0, 13, 38, 5', 24', 32'$], [$0, 7', 8', 13', 17', 38'$]]
        \item \{0=1060, 1=42824, 2=514842, 3=1945100, 4=1948174\}\newline [[$0, 1, 3, 15, 37, 0'$], [$0, 4, 24, 47, 22', 48'$], [$0, 5, 45, 9', 15', 25'$], [$0, 7, 28, 34', 36', 39'$],\newline  [$0, 9, 27, 14', 21', 46'$], [$0, 11, 3', 7', 41', 42'$], [$0, 2', 13', 26', 35', 43'$]]
        \item \{0=1590, 1=43672, 2=551412, 3=1934712, 4=1920614\}\newline [[$0, 1, 3, 17, 49, 0'$], [$0, 6, 29, 40, 48', 49'$], [$0, 8, 33, 15', 18', 34'$], [$0, 9, 27, 2', 30', 40'$],\newline  [$0, 10, 22, 33', 47', 51'$], [$0, 15, 6', 27', 32', 39'$], [$0, 5', 14', 16', 22', 45'$]]\item \{0=1590, 1=34344, 2=493854, 3=1970964, 4=1951248\}\newline [[$0, 1, 3, 22, 48, 0'$], [$0, 4, 24, 42, 10', 38'$], [$0, 7, 43, 15', 19', 33'$], [$0, 9, 23, 11', 32', 45'$],\newline  [$0, 12, 37, 1', 28', 30'$], [$0, 13, 7', 37', 40', 48'$], [$0, 3', 4', 13', 20', 51'$]]
        \item \{0=1590, 1=38584, 2=545688, 3=1943616, 4=1922522\}\newline [[$0, 1, 3, 10, 21, 0'$], [$0, 4, 17, 15', 42', 45'$], [$0, 5, 28, 34, 1', 10'$], [$0, 8, 22, 3', 16', 44'$],\newline  [$0, 12, 27, 31', 33', 39'$], [$0, 16, 9', 23', 30', 40'$], [$0, 2', 13', 17', 18', 37'$]]
        \item \{0=1590, 1=28408, 2=499578, 3=1960788, 4=1961636\}\newline [[$0, 1, 3, 11, 31, 0'$], [$0, 4, 41, 7', 24', 29'$], [$0, 5, 14, 32, 37', 45'$], [$0, 6, 40, 1', 4', 8'$],\newline  [$0, 7, 36, 16', 18', 46'$], [$0, 15, 6', 30', 43', 49'$], [$0, 12', 26', 27', 38', 47'$]]
        \item \{0=2120, 1=40280, 2=556182, 3=1934500, 4=1918918\}\newline [[$0, 1, 3, 9, 22, 0'$], [$0, 4, 15, 33, 20', 36'$], [$0, 5, 12, 28, 1', 23'$], [$0, 10, 22', 25', 27', 39'$],\newline  [$0, 14, 4', 8', 38', 51'$], [$0, 17, 10', 19', 30', 45'$], [$0, 26, 6', 7', 14', 35'$]]
    \end{enumerate}
    \end{example}

    \section{Steiner systems S(2,6,111)}

    In this section we present the Steiner systems $S(2,6,111)$ generated by an action of the cyclic group $\mathbb Z_{55}$ with two orbits of length $55$ and one fixed point, denoted by $\infty$. Elements of the first orbit are identified with the elements $0,\dots,54$ of the group $\mathbb Z_{55}$; elements of the second orbit are enumerated by primed numbers $0',\dots,54'$. Each designs is preceded by its fingerprint.

    \begin{example} There exist exactly four designs for $v=111$ generated by the action of the cyclic group $\mathbb Z_{55}$ with two orbits of size $55$ and one fixed point.
    \begin{enumerate}
        \item \{1=32120, 2=528000, 3=2207040, 4=2361040\}\quad [[$0, 11, 22, 33, 44, \infty$], [$0', 11', 22', 33', 44', \infty$], [$0, 1, 4, 49, 0', 13'$], [$0, 2, 15, 27, 31', 37'$], [$0, 5, 19, 39, 30', 47'$], [$0, 8, 31, 32', 34', 48'$],\newline [$0, 9, 38, 27', 45', 52'$], [$0, 18, 2', 23', 33', 38'$], [$0, 21', 41', 49', 50', 53'$]]
        \item \{1=33880, 2=534270, 3=2195820, 4=2364230\}\quad [[$0, 11, 22, 33, 44, \infty$], [$0', 11', 22', 33', 44', \infty$], [$0, 1, 3, 17, 0', 26'$], [$0, 4, 12, 2', 17', 22'$], [$0, 5, 29, 11', 41', 48'$], [$0, 6, 27, 46, 7', 30'$],\newline [$0, 7, 32, 42, 21', 27'$], [$0, 18, 33', 46', 47', 49'$], [$0, 4', 8', 32', 42', 51'$]]
        \item \{1=36960, 2=539550, 3=2193180, 4=2358510\}\quad [[$0, 11, 22, 33, 44, \infty$], [$0', 11', 22', 33', 44', \infty$], [$0, 1, 3, 32, 0', 12'$], [$0, 4, 10, 19, 44', 46'$], [$0, 5, 39, 1', 8', 21'$], [$0, 7, 20', 29', 45', 48'$],\newline [$0, 8, 20, 38, 14', 15'$], [$0, 13, 41, 5', 39', 43'$], [$0, 4', 10', 18', 28', 33'$]]
        \item \{0=550, 1=34320, 2=507540, 3=2224640, 4=2361150\}\quad [[$0, 11, 22, 33, 44, \infty$],\newline [$0', 11', 22', 33', 44', \infty$], [$0, 1, 3, 46, 0', 15'$], [$0, 4, 35, 41, 11', 21'$], [$0, 5, 28, 1', 18', 44'$],\newline [$0, 7, 26, 47, 29', 53'$], [$0, 13, 30, 5', 23', 32'$], [$0, 16, 36', 42', 49', 50'$], [$0, 4', 8', 38', 40', 43'$]]
    \end{enumerate}
    \end{example}

    Next construction is very rich in design generation. This designs are generated by an action of the group $G=\IZ_{37}\rtimes \IZ_3$ with three orbits of size $37,37,37$. The points of orbits are enumerated by numbers $0_i,\dots,36_i$ where $i$ denotes the number of an orbit. The group operation is defined by the formula $(a,b)\cdot (x,y)=(a+10^b\cdot x\mod 37, b+y \mod 3)$. The action of the group $G$ on the orbits are determined by the formula $(a,b)\cdot x=a+10^b\cdot x\mod 37$, where $(a,b)\in G=\IZ_{37}\times \IZ_3$ and $x\in \IZ_{37}$.

    The total number of non-isomorphic designs for construction above is $513$. In this paper we present ten of them, the other designs can be accessed by the URL \url{https://github.com/Ihromant/math-utils/blob/master/src/test/resources/ref/6-111-37%2C37%2C37.txt}

    \begin{example} Ten of $513$ designs $S(2,6,111)$ with three orbits of size $37,37,37$.
    \begin{enumerate}
        \item \{1=27528, 2=504828, 3=2212008, 4=2383836\} \quad [$[0_0, 1_0, 4_0, 11_0, 4_1, 4_2]$,\newline $[0_0, 2_0, 7_1, 11_1, 24_1, 12_2]$, $[0_0, 5_0, 24_0, 15_1, 25_1, 26_1]$, $[0_0, 6_0, 14_1, 11_2, 20_2, 28_2]$, \newline $[0_0, 9_0, 25_0, 6_2, 32_2, 33_2]$, $[0_1, 6_1, 29_1, 4_2, 9_2, 22_2]$, $[0_1, 9_1, 25_1, 20_2, 24_2, 27_2]$]
        \item \{1=28416, 2=546786, 3=2180484, 4=2372514\} \quad [$[0_0, 1_0, 4_0, 11_0, 4_1, 4_2]$,\newline $[0_0, 2_0, 14_1, 17_1, 28_1, 34_2]$, $[0_0, 5_0, 24_0, 6_2, 10_2, 13_2]$, $[0_0, 6_0, 14_0, 11_1, 13_1, 33_1]$, \newline $[0_0, 9_0, 32_1, 20_2, 25_2, 31_2]$, $[0_1, 5_1, 18_1, 10_2, 19_2, 31_2]$, $[0_1, 9_1, 21_1, 16_2, 18_2, 33_2]$]
        \item \{1=31080, 2=554778, 3=2175156, 4=2367186\}\quad  [$[0_0, 1_0, 4_0, 11_0, 4_1, 4_2]$,\newline $[0_0, 2_0, 22_0, 9_1, 18_1, 34_1]$, $[0_0, 5_0, 18_0, 29_2, 32_2, 36_2]$, $[0_0, 6_0, 27_1, 8_2, 16_2, 34_2]$,\newline  $[0_0, 9_0, 25_0, 1_1, 14_1, 19_1]$, $[0_0, 2_1, 6_1, 17_1, 5_2, 22_2]$, $[0_1, 6_1, 29_1, 14_2, 23_2, 35_2]$]
        \item \{1=34632, 2=534132, 3=2179152, 4=2380284\}\quad [$[0_0, 1_0, 4_0, 11_0, 4_1, 4_2]$,\newline $[0_0, 2_0, 17_0, 12_1, 18_1, 26_1]$, $[0_0, 5_0, 24_0, 14_2, 16_2, 36_2]$, $[0_0, 6_0, 25_1, 23_2, 26_2, 34_2]$, \newline $[0_0, 9_0, 25_0, 8_1, 11_1, 15_1]$, $[0_0, 7_1, 22_1, 31_1, 13_2, 18_2]$, $[0_1, 1_1, 11_1, 3_2, 15_2, 31_2]$]
        \item \{1=37296, 2=576756, 3=2180040, 4=2334108\} \quad [$[0_0, 1_0, 11_0, 6_1, 19_1, 24_1]$,\newline $[0_0, 2_0, 17_0, 15_2, 19_2, 22_2]$, $[0_0, 3_0, 12_0, 33_0, 12_1, 12_2]$, $[0_0, 5_0, 24_0, 1_1, 31_1, 34_1]$, \newline $[0_0, 6_0, 28_1, 3_2, 27_2, 29_2]$, $[0_0, 2_1, 3_1, 11_1, 26_2, 32_2]$, $[0_1, 2_1, 22_1, 7_2, 19_2, 35_2]$]
        \item \{1=38184, 2=509490, 3=2218668, 4=2361858\} \quad [$[0_0, 1_0, 11_0, 3_1, 15_1, 31_1]$,\newline $[0_0, 2_0, 17_0, 12_1, 18_1, 26_1]$, $[0_0, 3_0, 11_1, 1_2, 15_2, 36_2]$, $[0_0, 5_0, 14_0, 33_1, 18_2, 28_2]$, \newline $[0_0, 0_1, 17_1, 22_1, 35_1, 0_2]$, $[0_1, 1_1, 27_1, 6_2, 9_2, 13_2]$, $[0_1, 3_1, 33_1, 10_2, 29_2, 34_2]$]
        \item \{1=44400, 2=564102, 3=2181372, 4=2338326\}\quad [$[0_0, 1_0, 4_0, 11_0, 4_1, 4_2]$,\newline $[0_0, 2_0, 22_0, 10_1, 23_1, 28_1]$, $[0_0, 5_0, 16_1, 1_2, 25_2, 29_2]$, $[0_0, 6_0, 14_0, 12_2, 22_2, 23_2]$, \newline $[0_0, 9_0, 25_0, 5_1, 7_1, 22_1]$, $[0_0, 2_1, 31_1, 32_1, 5_2, 27_2]$, $[0_1, 9_1, 21_1, 15_2, 23_2, 29_2]$]
        \item \{0=1110, 1=43512, 2=555444, 3=2197800, 4=2330334\}\; [$[0_0, 1_0, 4_0, 11_0, 4_1, 4_2]$,\newline $[0_0, 2_0, 17_0, 7_1, 13_1, 36_1]$, $[0_0, 5_0, 25_1, 13_2, 27_2, 29_2]$, $[0_0, 6_0, 14_0, 15_2, 16_2, 26_2]$, \newline $[0_0, 9_0, 1_1, 17_1, 18_1, 34_2]$, $[0_1, 3_1, 33_1, 14_2, 27_2, 32_2]$, $[0_1, 5_1, 24_1, 6_2, 10_2, 13_2]$]
        \item \{0=1110, 1=44400, 2=538794, 3=2199132, 4=2344764\}\; [$[0_0, 1_0, 11_0, 6_1, 19_1, 24_1]$,\newline $[0_0, 2_0, 17_0, 1_1, 27_1, 28_1]$, $[0_0, 3_0, 24_0, 17_1, 11_2, 33_2]$, $[0_0, 6_0, 29_0, 20_2, 21_2, 31_2]$, \newline $[0_0, 0_1, 9_1, 12_1, 16_1, 0_2]$, $[0_0, 2_1, 33_1, 1_2, 19_2, 35_2]$, $[0_1, 2_1, 22_1, 3_2, 26_2, 32_2]$]
        \item \{0=1110, 1=44400, 2=572760, 3=2208456, 4=2301474\}\;  [$[0_0, 1_0, 4_0, 11_0, 4_1, 4_2]$,\newline $[0_0, 2_0, 17_0, 12_1, 18_1, 26_1]$, $[0_0, 5_0, 34_1, 21_2, 24_2, 36_2]$, $[0_0, 6_0, 5_1, 21_1, 23_1, 26_2]$, \newline $[0_0, 9_0, 21_0, 6_2, 7_2, 17_2]$, $[0_1, 1_1, 11_1, 9_2, 17_2, 23_2]$, $[0_1, 3_1, 33_1, 10_2, 29_2, 34_2]$]
    \end{enumerate}
    \end{example}

    \section{Steiner systems S(2,6,96)}

    In this section we present all 133 Steiner designs $S(2,6,96)$ generated by the action of the group $G=\IZ_{19}\rtimes \IZ_3$ with four orbits of size $57,19,19,1$. The points of the first orbit are enumerated by numbers $0,\dots,56$; the points of the second and third orbits are enumerated by numbers $0',\dots,18'$ and $0'',\dots,18''$ with one or two primes; the fourth orbit contains a unique fixed point denoted by $\infty$.

    Elements of the group $G$ are enumerated with numbers $0,\dots,56$. Each number $n\in\{0,\dots,56\}=G$ is identified with the unique pair $(x,y)\in\IZ_{19}\times\IZ_3$ such that $n=3x+y$. The group operation is defined by the formula $(a,b)\cdot (x,y)=(a+7^b\cdot x\mod 19, b+y \mod 3)$, which also determines the action of the group $G$ on points of the first orbit. The action of the group $G$ on the second and the third orbits are determined by the formula $(a,b)\cdot x=a+7^b\cdot x\mod 19$, where $(a,b)\in G=\IZ_{19}\times \IZ_3$ and $x\in \IZ_{19}$.

    A simpler construction, but richer in results appears for the same group $G=\IZ_{19}\rtimes \IZ_3$ acting on a set with six orbits of cardinality   $19,19,19,19,19,1$. The points of five $19$-element orbits are denoted by the numbers $0_i,\dots,18_i$ indexed by the number $i$ of an orbit. The action of the group $G$ on the $19$-element orbits is defined in the preceding paragraph. The last orbit contains a unique fixed point denoted by $\infty$. The total number of discovered Steiner designs $S(2,6,96)$ equals $634$; all of them can be accessed by the URL \url{https://github.com/Ihromant/math-utils/blob/master/src/test/resources/ref/6-96-1%2C19%2C19%2C19%2C19%2C19.txt}

    \begin{example} Ten of $634$ designs $S(2,6,96)$ with six orbits of size $19,19,19,19,19,1$.
    \begin{enumerate}
        \item \{1=43776, 2=467514, 3=1519620, 4=1252290\}\;
        [$[0_0,\! 1_0,\! 8_0,\! 4_1,\! 10_1,\! 14_1]$, $[0_0,\! 2_0,\! 16_0,\! 7_4,\! 13_4,\! 17_4]$,\newline $[0_0, 4_0, 10_0, 8_3, 9_3, 16_3]$, $[0_0, 0_1, 0_2, 0_3, 0_4, \infty]$, $[0_0, 1_1, 12_2, 13_2, 7_3, 3_4]$, $[0_0, 5_1, 2_2, 7_2, 13_3, 18_4]$, $[0_0, 8_1, 4_2, 17_2, 3_3, 6_4]$, $[0_1, 1_1, 8_1, 7_4, 9_4, 12_4]$, $[0_1, 2_1, 16_1, 7_3, 13_3, 17_3]$, $[0_3, 2_3, 16_3, 1_4, 8_4, 9_4]$]
        \item \{1=33288, 2=456228, 3=1507536, 4=1286148\}\; [$[0_0, 1_0, 8_0, 13_1, 16_1, 18_1]$,\newline $[0_0, 2_0, 1_3, 17_3, 12_4, 16_4]$, $[0_0, 4_0, 10_0, 8_2, 9_2, 16_2]$, $[0_0, 0_1, 0_2, 0_3, 0_4, \infty]$, $[0_0, 1_1, 14_1, 10_2, 9_4, 11_4]$, $[0_0, 4_1, 1_2, 3_2, 6_3, 14_3]$, $[0_1, 1_1, 12_1, 5_3, 9_3, 18_3]$, $[0_1, 1_2, 14_2, 7_3, 4_4, 11_4]$]
        \item \{1=33288, 2=446994, 3=1506396, 4=1296522\}\; [$[0_0,\! 1_0,\! 8_0,\! 2_1,\! 11_1,\! 15_1]$, $[0_0,\! 2_0,\! 16_0,\! 1_4,\! 8_4,\! 9_4]$,\newline $[0_0, 4_0, 10_0, 6_3, 13_3, 14_3]$, $[0_0, 0_1, 0_2, 0_3, 0_4, \infty]$, $[0_0, 4_1, 9_2, 13_2, 16_3, 3_4]$, $[0_0, 5_1, 8_2, 16_2, 1_3, 15_4]$, $[0_0, 8_1, 2_2, 7_2, 12_3, 5_4]$, $[0_1, 1_1, 12_1, 2_4, 4_4, 7_4]$, $[0_1, 2_1, 5_1, 3_3, 7_3, 16_3]$, $[0_3, 2_3, 16_3, 7_4, 13_4, 17_4]$]
        \item \{0=570, 1=41040, 2=472986, 3=1517796, 4=1250808\}\newline [$[0_0, 1_0, 8_0, 2_1, 11_1, 15_1]$, $[0_0, 2_0, 16_0, 5_4, 14_4, 18_4]$, $[0_0, 4_0, 13_0, 3_2, 6_2, 8_2]$, $[0_0, 0_1, 0_2, 0_3, 0_4, \infty]$, $[0_0, 4_1, 12_1, 1_2, 5_3, 6_3]$, $[0_0, 5_1, 2_3, 11_3, 1_4, 13_4]$, $[0_0, 5_2, 13_2, 12_3, 15_3, 9_4]$, $[0_1, 2_1, 4_2, 13_2, 11_4, 16_4]$]
        \item \{1=34656, 2=458964, 3=1515744, 4=1273836\}\newline [$[0_0, 1_0, 8_0, 4_1, 10_1, 14_1]$, $[0_0, 2_0, 16_0, 3_3, 4_3, 11_3]$, $[0_0, 4_0, 1_1, 15_2, 16_2, 3_4]$, $[0_0, 0_1, 0_2, 0_3, 0_4, \infty]$, $[0_0, 8_1, 5_3, 6_4, 15_4, 17_4]$, $[0_0, 2_2, 4_2, 15_3, 18_3, 1_4]$, $[0_1, 1_1, 8_1, 2_3, 11_3, 15_3]$, $[0_1, 2_1, 9_2, 18_2, 8_3, 12_4]$]
        \item \{0=1140, 1=36936, 2=472302, 3=1509132, 4=1263690\}\newline [$[0_0, 1_0, 8_0, 2_1, 11_1, 15_1]$, $[0_0, 2_0, 16_0, 1_4, 8_4, 9_4]$, $[0_0, 4_0, 10_0, 8_3, 9_3, 16_3]$, $[0_0, 0_1, 0_2, 0_3, 0_4, \infty]$, $[0_0, 4_1, 6_2, 11_2, 1_3, 5_4]$, $[0_0, 5_1, 3_2, 15_2, 14_3, 10_4]$, $[0_0, 8_1, 12_2, 16_2, 10_3, 2_4]$,\newline $[0_1, 1_1, 12_1, 8_3, 11_3, 13_3]$, $[0_1, 2_1, 5_1, 4_4, 8_4, 14_4]$, $[0_3, 4_3, 13_3, 2_4, 16_4, 18_4]$]
        \item \{0=1140, 1=43776, 2=487692, 3=1489752, 4=1260840\}\newline [$[0_0, 1_0, 8_0, 4_1, 10_1, 14_1]$, $[0_0, 2_0, 16_0, 3_3, 4_3, 11_3]$, $[0_0, 4_0, 5_2, 17_3, 2_4, 15_4]$, $[0_0, 0_1, 0_2, 0_3, 0_4, \infty]$, $[0_0, 1_1, 18_1, 3_2, 4_4, 12_4]$, $[0_0, 5_1, 4_2, 10_2, 12_2, 8_3]$, $[0_1, 1_1, 8_1, 7_3, 9_3, 12_3]$, $[0_1, 10_2, 13_3, 17_3, 9_4, 12_4]$]
        \item \{0=1140, 1=46512, 2=482562, 3=1500924, 4=1252062\}\newline [$[0_0, 1_0, 8_0, 2_1, 11_1, 15_1]$, $[0_0, 2_0, 16_0, 1_3, 8_3, 9_3]$, $[0_0, 4_0, 10_0, 6_2, 13_2, 14_2]$, $[0_0, 0_1, 0_2, 0_3, 0_4, \infty]$, $[0_0, 4_1, 12_1, 17_3, 3_4, 10_4]$, $[0_0, 5_1, 11_2, 17_2, 7_4, 12_4]$, $[0_0, 8_2, 10_3, 14_3, 6_4, 16_4]$, $[0_1, 2_1, 3_2, 17_2, 11_3, 14_3]$]
        \item \{0=3990, 1=49704, 2=478458, 3=1497276, 4=1253772\}\newline [$[0_0, 1_0, 8_0, 2_1, 11_1, 15_1]$, $[0_0, 2_0, 16_0, 3_3, 4_3, 11_3]$, $[0_0, 4_0, 13_0, 9_4, 10_4, 17_4]$, $[0_0, 0_1, 0_2, 0_3, 0_4, \infty]$, $[0_0, 4_1, 18_1, 3_2, 9_2, 12_3]$, $[0_0, 5_1, 8_2, 16_2, 3_4, 7_4]$, $[0_0, 1_2, 15_2, 13_3, 16_3, 12_4]$, $[0_1, 1_1, 2_3, 6_3, 9_4, 11_4]$]
        \item \{0=570, 1=41952, 2=457596, 3=1512096, 4=1270986\}\newline [$[0_0, 1_0, 8_0, 2_1, 11_1, 15_1]$, $[0_0, 2_0, 16_0, 1_4, 8_4, 9_4]$, $[0_0, 4_0, 10_0, 1_3, 15_3, 17_3]$, $[0_0, 0_1, 0_2, 0_3, 0_4, \infty]$, $[0_0, 4_1, 12_2, 17_2, 18_3, 3_4]$, $[0_0, 5_1, 2_2, 9_2, 14_3, 10_4]$, $[0_0, 8_1, 11_2, 15_2, 6_3, 17_4]$, $[0_1, 1_1, 8_1, 2_4, 11_4, 15_4]$, $[0_1, 2_1, 5_1, 1_3, 12_3, 13_3]$, $[0_3, 4_3, 13_3, 2_4, 16_4, 18_4]$]
    \end{enumerate}
    \end{example}

    \begin{example} Designs with four orbits of size $57,19,19,1$.
    \begin{enumerate}
        \item \{1=25992, 2=432630, 3=1535580, 4=1288998\} \newline [[$0, 1, 2, 0', 0'', \infty$], [$0, 8, 16, 19, 29, 33$], [$3, 7, 31, 38, 47, 0'$], [$5, 12, 35, 50, 0'', 2''$], \newline [$14, 43, 0', 2', 11'', 18''$], [$15, 20, 52, 0', 1', 8'$], [$16, 25, 46, 0', 4', 14''$], [$26, 27, 49, 0'', 4'', 10''$]]
        \item \{1=28728, 2=451440, 3=1521672, 4=1281360\} \newline [[$0, 1, 2, 0', 0'', \infty$], [$0, 6, 19, 34, 38, 47$], [$3, 34, 39, 48, 0', 2'$], [$6, 50, 55, 0', 4', 13'$], \newline [$9, 29, 36, 0', 2'', 6''$], [$15, 37, 39, 44, 0'', 1''$], [$27, 30, 45, 53, 0', 1''$], [$0', 1', 8', 13'', 16'', 18''$]]
        \item \{1=28728, 2=452124, 3=1529880, 4=1272468\} \newline [[$0, 1, 2, 0', 0'', \infty$], [$0, 3, 11, 19, 35, 55$], [$4, 12, 29, 48, 55, 0'$], [$6, 35, 46, 0', 1', 8'$], \newline [$8, 17, 47, 53, 0'', 1''$], [$9, 31, 0', 4', 1'', 14''$], [$16, 42, 56, 0', 2', 16'$], [$30, 34, 53, 0', 4'', 18''$]]
        \item \{1=31464, 2=449730, 3=1527372, 4=1274634\} \newline [[$0, 1, 2, 0', 0'', \infty$], [$3, 27, 37, 43, 0', 4''$], [$4, 9, 14, 16, 40, 0''$], [$11, 12, 34, 0'', 2'', 16''$], \newline [$13, 32, 42, 0'', 1'', 8''$], [$15, 30, 34, 42, 0', 2'$], [$29, 45, 51, 0', 4', 16''$], [$0', 1', 8', 2'', 11'', 15''$]]
        \item \{1=31920, 2=452124, 3=1539456, 4=1259700\} \newline [[$0, 1, 2, 0', 0'', \infty$], [$0, 4, 12, 26, 31, 44$], [$0, 5, 7, 25, 27, 50$], [$5, 42, 0', 4', 15'', 18''$], \newline [$6, 35, 46, 0', 1', 8'$], [$15, 18, 31, 54, 0', 2'$], [$27, 33, 55, 0', 4'', 5''$], [$35, 39, 50, 56, 0'', 4''$]]
        \item \{1=32832, 2=457254, 3=1526460, 4=1266654\} \newline [[$0, 1, 2, 0', 0'', \infty$], [$0, 4, 8, 16, 38, 39$], [$7, 25, 50, 54, 0'', 4''$], [$8, 14, 36, 47, 50, 0'$], \newline [$12, 31, 44, 0', 1', 8'$], [$18, 38, 43, 0', 1'', 18''$], [$23, 39, 52, 0', 2', 16'$], [$34, 49, 0', 4', 2'', 10''$]]
        \item \{1=34656, 2=456570, 3=1512324, 4=1279650\} \newline [[$0, 1, 2, 0', 0'', \infty$], [$0, 11, 22, 34, 38, 39$], [$10, 20, 24, 0', 4', 13'$], [$13, 31, 39, 55, 0'', 1''$], \newline [$15, 22, 25, 47, 0', 2'$], [$17, 37, 43, 0', 1'', 17''$], [$18, 21, 33, 46, 0', 8''$], [$0', 1', 8', 4'', 10'', 14''$]]
        \item \{1=34656, 2=478116, 3=1494768, 4=1275660\} \newline [[$0, 1, 2, 0', 0'', \infty$], [$0, 11, 22, 34, 38, 39$], [$7, 14, 20, 32, 49, 0'$], [$10, 16, 47, 0', 4', 18''$], \newline [$12, 15, 40, 48, 0'', 2''$], [$17, 31, 0', 2', 11'', 17''$], [$21, 29, 37, 0', 1', 8'$], [$25, 35, 39, 0'', 1'', 12''$]]
        \item \{1=35112, 2=462384, 3=1528512, 4=1257192\} \newline [[$0, 1, 2, 0', 0'', \infty$], [$3, 36, 39, 53, 54, 0''$], [$6, 26, 36, 55, 0', 5''$], [$10, 14, 33, 0'', 2'', 16''$], \newline [$20, 40, 42, 0', 4', 10'$], [$27, 35, 39, 50, 0', 2'$], [$31, 48, 53, 0', 1'', 8''$], [$0', 1', 8', 4'', 10'', 14''$]]
        \item \{1=35568, 2=444942, 3=1505940, 4=1296750\} \newline [[$0, 1, 2, 0', 0'', \infty$], [$0, 3, 14, 34, 38, 55$], [$6, 19, 28, 31, 50, 0'$], [$7, 14, 21, 53, 0'', 1''$], \newline [$9, 15, 36, 0', 4'', 8''$], [$21, 29, 37, 0', 1', 8'$], [$22, 42, 47, 0', 4', 13'$], [$23, 38, 0', 2', 13'', 16''$]]
        \item \{1=35568, 2=458964, 3=1512096, 4=1276572\} \newline [[$0, 1, 2, 0', 0'', \infty$], [$0, 4, 13, 15, 35, 41$], [$4, 19, 24, 42, 51, 0'$], [$6, 35, 46, 0', 1', 8'$], \newline [$9, 12, 56, 0', 5'', 12''$], [$14, 17, 21, 35, 0'', 2''$], [$21, 37, 0', 2', 4'', 13''$], [$29, 30, 52, 0', 4', 13'$]]
        \item \{1=36024, 2=446310, 3=1528740, 4=1272126\} \newline [[$0, 1, 2, 0', 0'', \infty$], [$0, 4, 26, 41, 48, 55$], [$8, 19, 35, 0', 4', 14''$], [$12, 24, 40, 54, 0', 2'$], \newline [$12, 31, 44, 0'', 1'', 8''$], [$15, 18, 38, 39, 0'', 4''$], [$25, 32, 45, 51, 0', 5''$], [$0', 1', 8', 7'', 9'', 12''$]]
        \item \{1=36480, 2=456228, 3=1520760, 4=1269732\} \newline [[$0, 1, 2, 0', 0'', \infty$], [$0, 4, 8, 26, 34, 45$], [$3, 14, 44, 47, 0'', 2''$], [$13, 19, 28, 45, 0', 8''$], \newline [$15, 22, 32, 42, 0', 2'$], [$16, 18, 53, 0'', 1'', 8''$], [$33, 52, 55, 0', 4', 9''$], [$0', 1', 8', 2'', 11'', 15''$]]
        \item \{1=36936, 2=449388, 3=1519392, 4=1277484\} \newline [[$0, 1, 2, 0', 0'', \infty$], [$0, 10, 17, 19, 23, 42$], [$3, 6, 14, 30, 56, 0'$], [$5, 10, 0', 4', 7'', 13''$], \newline [$13, 24, 44, 0'', 2'', 5''$], [$13, 37, 39, 0', 2', 18''$], [$14, 25, 37, 52, 0'', 1''$], [$15, 20, 52, 0', 1', 8'$]]
        \item \{1=37392, 2=461358, 3=1519620, 4=1264830\} \newline [[$0, 1, 2, 0', 0'', \infty$], [$4, 6, 9, 37, 51, 0''$], [$8, 37, 45, 50, 0', 4'$], [$9, 13, 35, 0', 2', 16'$], \newline [$11, 20, 39, 47, 0', 4''$], [$27, 34, 52, 0', 1'', 2''$], [$29, 30, 52, 0'', 4'', 13''$], [$0', 1', 8', 13'', 16'', 18''$]]
        \item \{1=37392, 2=484272, 3=1504800, 4=1256736\} \newline [[$0, 1, 2, 0', 0'', \infty$], [$3, 6, 17, 31, 42, 0''$], [$9, 19, 26, 0'', 4'', 13''$], [$15, 19, 37, 39, 0', 2'$], \newline [$15, 28, 56, 0'', 1'', 12''$], [$22, 36, 52, 0', 4', 2''$], [$28, 34, 43, 51, 0', 10''$], [$0', 1', 8', 7'', 9'', 12''$]]
        \item \{1=37392, 2=490428, 3=1503888, 4=1251492\} \newline [[$0, 1, 2, 0', 0'', \infty$], [$0, 5, 22, 28, 47, 48$], [$4, 12, 37, 47, 53, 0''$], [$5, 17, 38, 0', 2'', 8''$], \newline [$6, 35, 46, 0', 1', 8'$], [$9, 13, 35, 0'', 2'', 16''$], [$16, 20, 29, 47, 0', 2'$], [$49, 56, 0', 4', 5'', 13''$]]
        \item \{1=38304, 2=468882, 3=1501380, 4=1274634\} \newline [[$0, 1, 2, 0', 0'', \infty$], [$0, 4, 12, 26, 31, 44$], [$7, 12, 40, 46, 49, 0'$], [$10, 22, 24, 53, 0'', 1''$], \newline [$21, 29, 37, 0', 1', 8'$], [$22, 43, 50, 0', 4', 6''$], [$23, 25, 0', 2', 13'', 18''$], [$29, 30, 52, 0'', 4'', 13''$]]
        \item \{1=38304, 2=495558, 3=1509132, 4=1240206\} \newline [[$0, 1, 2, 0', 0'', \infty$], [$0, 20, 22, 25, 33, 56$], [$4, 25, 31, 49, 51, 0'$], [$6, 35, 46, 0', 1', 8'$], \newline [$9, 52, 0', 4', 5'', 6''$], [$10, 15, 19, 23, 0'', 4''$], [$13, 42, 55, 0', 8'', 10''$], [$20, 27, 34, 0', 2', 5'$]]
        \item \{1=38304, 2=498636, 3=1489296, 4=1256964\} \newline [[$0, 1, 2, 0', 0'', \infty$], [$0, 6, 35, 40, 44, 55$], [$0, 7, 14, 17, 21, 31$], [$9, 33, 51, 0', 1'', 16''$], \newline [$11, 34, 0', 4', 10'', 18''$], [$14, 24, 32, 47, 0', 2'$], [$15, 20, 52, 0', 1', 8'$], [$15, 34, 45, 56, 0'', 2''$]]
        \item \{1=38760, 2=449730, 3=1508676, 4=1286034\} \newline [[$0, 1, 2, 0', 0'', \infty$], [$5, 6, 14, 32, 48, 0''$], [$5, 24, 43, 51, 0', 2''$], [$15, 18, 25, 39, 0', 2'$], \newline [$20, 40, 42, 0'', 4'', 10''$], [$22, 42, 47, 0', 4', 13'$], [$31, 37, 48, 0', 5'', 13''$], [$0', 1', 8', 7'', 9'', 12''$]]
        \item \{1=38760, 2=470934, 3=1510044, 4=1263462\} \newline [[$0, 1, 2, 0', 0'', \infty$], [$0, 6, 13, 28, 47, 56$], [$3, 6, 37, 46, 0', 4'$], [$5, 18, 39, 51, 55, 0''$], \newline [$12, 31, 44, 0', 1', 8'$], [$15, 33, 42, 0', 2'', 9''$], [$19, 30, 0', 2', 13'', 18''$], [$20, 40, 42, 0'', 4'', 10''$]]
        \item \{1=38760, 2=472644, 3=1506168, 4=1265628\} \newline [[$0, 1, 2, 0', 0'', \infty$], [$0, 4, 14, 22, 46, 52$], [$5, 25, 29, 0', 4'', 17''$], [$6, 35, 46, 0', 1', 8'$], \newline [$16, 23, 44, 0', 2', 12''$], [$23, 28, 36, 0'', 1'', 8''$], [$26, 40, 53, 0', 4', 7''$], [$30, 53, 55, 0'', 2'', 5''$]]
        \item \{1=38760, 2=474354, 3=1511412, 4=1258674\} \newline [[$0, 1, 2, 0', 0'', \infty$], [$3, 8, 20, 34, 56, 0''$], [$4, 52, 0', 2', 4'', 18''$], [$6, 34, 42, 0', 4', 15''$], \newline [$12, 31, 44, 0', 1', 8'$], [$13, 19, 36, 49, 54, 0'$], [$17, 19, 51, 0'', 1'', 8''$], [$23, 43, 45, 0'', 4'', 13''$]]
        \item \{1=38760, 2=474696, 3=1512552, 4=1257192\} \newline [[$0, 1, 2, 0', 0'', \infty$], [$0, 4, 21, 26, 29, 37$], [$3, 10, 31, 36, 54, 0'$], [$4, 13, 43, 51, 0', 2''$], \newline [$20, 40, 46, 0', 4'', 18''$], [$25, 35, 39, 0'', 1'', 12''$], [$27, 52, 0', 1', 6', 11''$], [$29, 30, 52, 0'', 4'', 13''$]]
        \item \{1=38760, 2=481878, 3=1513692, 4=1248870\} \newline [[$0, 1, 2, 0', 0'', \infty$], [$0, 9, 23, 34, 38, 40$], [$7, 17, 21, 0', 4', 10'$], [$10, 22, 43, 51, 0', 5''$], \newline [$18, 33, 37, 53, 0', 2'$], [$24, 36, 46, 0', 1'', 18''$], [$29, 38, 45, 50, 0'', 1''$], [$0', 1', 8', 4'', 10'', 14''$]]
        \item \{1=38760, 2=507186, 3=1486788, 4=1250466\} \newline [[$0, 1, 2, 0', 0'', \infty$], [$0, 8, 16, 18, 28, 35$], [$3, 15, 47, 0', 2', 9''$], [$6, 27, 32, 47, 0'', 1''$], \newline [$8, 49, 54, 0'', 4'', 13''$], [$17, 18, 24, 51, 0', 4'$], [$23, 27, 42, 56, 0', 2''$], [$0', 1', 8', 13'', 16'', 18''$]]
        \item \{1=39216, 2=468198, 3=1521444, 4=1254342\} \newline [[$0, 1, 2, 0', 0'', \infty$], [$4, 31, 36, 40, 42, 0'$], [$13, 21, 38, 42, 0'', 2''$], [$14, 21, 35, 43, 0', 8''$], \newline [$15, 24, 54, 0', 2', 7''$], [$20, 31, 48, 0'', 1'', 12''$], [$28, 32, 51, 0', 4', 13'$], [$0', 1', 8', 4'', 10'', 14''$]]
        \item \{1=40128, 2=479826, 3=1492716, 4=1270530\} \newline [[$0, 1, 2, 0', 0'', \infty$], [$0, 8, 9, 23, 25, 31$], [$3, 21, 55, 0', 4'', 12''$], [$4, 6, 18, 42, 48, 0'$], \newline [$4, 15, 35, 0'', 4'', 13''$], [$12, 49, 0', 1', 6', 11''$], [$19, 26, 36, 40, 0', 2''$], [$32, 37, 45, 0'', 2'', 16''$]]
        \item \{1=41040, 2=449388, 3=1521216, 4=1271556\} \newline [[$0, 1, 2, 0', 0'', \infty$], [$0, 5, 20, 31, 37, 48$], [$4, 38, 47, 51, 0'', 4''$], [$5, 25, 27, 0', 2', 16'$], \newline [$6, 17, 29, 47, 53, 0'$], [$8, 16, 39, 0', 1'', 18''$], [$12, 31, 44, 0', 1', 8'$], [$35, 38, 0', 4', 6'', 17''$]]
        \item \{1=41040, 2=461358, 3=1517340, 4=1263462\} \newline [[$0, 1, 2, 0', 0'', \infty$], [$3, 7, 17, 36, 42, 0''$], [$9, 14, 26, 55, 0', 2'$], [$13, 21, 30, 46, 0', 5''$], \newline [$16, 25, 52, 0', 4'', 12''$], [$20, 27, 34, 0'', 2'', 5''$], [$23, 43, 45, 0', 4', 13'$], [$0', 1', 8', 2'', 11'', 15''$]]
        \item \{1=42408, 2=470250, 3=1504572, 4=1265970\} \newline [[$0, 1, 2, 0', 0'', \infty$], [$0, 7, 13, 41, 48, 56$], [$4, 13, 31, 35, 0', 8''$], [$10, 22, 25, 30, 0'', 2''$], \newline [$11, 28, 49, 0', 4', 5''$], [$20, 31, 48, 0'', 1'', 12''$], [$27, 44, 46, 48, 0', 2'$], [$0', 1', 8', 4'', 10'', 14''$]]
        \item \{1=42864, 2=454860, 3=1507992, 4=1277484\} \newline [[$0, 1, 2, 0', 0'', \infty$], [$0, 5, 9, 37, 43, 47$], [$4, 6, 16, 37, 46, 0''$], [$7, 48, 56, 0'', 4'', 13''$], \newline [$10, 32, 0', 2', 4'', 12''$], [$11, 27, 34, 52, 0', 4'$], [$12, 31, 44, 0', 1', 8'$], [$13, 16, 24, 0', 1'', 17''$]]
        \item \{1=42864, 2=463410, 3=1522356, 4=1254570\} \newline [[$0, 1, 2, 0', 0'', \infty$], [$0, 7, 23, 39, 50, 52$], [$3, 32, 40, 51, 54, 0'$], [$6, 20, 21, 33, 0', 8''$], \newline [$8, 18, 44, 0', 5'', 15''$], [$10, 23, 48, 0'', 2'', 5''$], [$24, 34, 41, 0'', 1'', 12''$], [$27, 48, 0', 1', 4', 7''$]]
        \item \{1=43320, 2=470250, 3=1507308, 4=1262322\} \newline [[$0, 1, 2, 0', 0'', \infty$], [$0, 4, 21, 26, 29, 37$], [$5, 17, 38, 43, 47, 0'$], [$8, 9, 26, 51, 0', 10''$], \newline [$12, 19, 35, 0', 1'', 8''$], [$19, 29, 33, 0'', 2'', 5''$], [$20, 40, 42, 0'', 4'', 10''$], [$36, 49, 0', 1', 6', 4''$]]
        \item \{1=43320, 2=494190, 3=1501836, 4=1243854\} \newline [[$0, 1, 2, 0', 0'', \infty$], [$5, 18, 27, 42, 54, 0''$], [$6, 16, 23, 0'', 4'', 10''$], [$6, 35, 46, 0', 1', 8'$], \newline [$13, 26, 27, 53, 55, 0'$], [$21, 28, 0', 4', 2'', 10''$], [$23, 42, 47, 0', 2', 1''$], [$31, 36, 47, 0'', 2'', 16''$]]
        \item \{1=43776, 2=467172, 3=1506168, 4=1266084\} \newline [[$0, 1, 2, 0', 0'', \infty$], [$4, 10, 43, 46, 54, 0''$], [$4, 45, 53, 0', 4', 10'$], [$8, 30, 48, 0', 1'', 12''$], \newline [$9, 34, 43, 56, 0', 2'$], [$14, 18, 32, 55, 0', 5''$], [$15, 44, 55, 0'', 2'', 16''$], [$0', 1', 8', 4'', 10'', 14''$]]
        \item \{1=48336, 2=451782, 3=1524180, 4=1258902\} \newline [[$0, 1, 2, 0', 0'', \infty$], [$4, 22, 38, 49, 0', 10''$], [$6, 21, 30, 38, 0'', 4''$], [$7, 28, 30, 36, 50, 0'$], \newline [$18, 41, 43, 0', 4', 10'$], [$25, 35, 39, 0'', 1'', 12''$], [$27, 32, 48, 0', 2', 5''$], [$0', 1', 8', 7'', 9'', 12''$]]
        \item \{0=570, 1=26448, 2=470250, 3=1511412, 4=1274520\} \newline [[$0, 1, 2, 0', 0'', \infty$], [$4, 15, 35, 0'', 4'', 13''$], [$6, 16, 23, 0', 4', 10'$], [$7, 26, 30, 46, 51, 0''$], \newline [$10, 18, 22, 25, 0', 4''$], [$15, 27, 30, 52, 0', 2'$], [$20, 31, 37, 0', 1'', 2''$], [$0', 1', 8', 13'', 16'', 18''$]]
        \item \{0=570, 1=31920, 2=437418, 3=1530108, 4=1283184\} \newline [[$0, 1, 2, 0', 0'', \infty$], [$0, 4, 16, 18, 26, 53$], [$0, 5, 22, 31, 42, 56$], [$6, 15, 36, 55, 0'', 1''$], \newline [$14, 28, 43, 47, 0', 2'$], [$21, 33, 44, 0', 1'', 17''$], [$29, 50, 0', 4', 8'', 14''$], [$39, 49, 56, 0', 1', 8'$]]
        \item \{0=570, 1=31920, 2=464778, 3=1506852, 4=1279080\} \newline [[$0, 1, 2, 0', 0'', \infty$], [$0, 20, 22, 25, 33, 56$], [$3, 17, 34, 0'', 4'', 13''$], [$4, 16, 18, 37, 46, 0''$], \newline [$14, 26, 55, 0', 2'', 5''$], [$15, 20, 52, 0', 1', 8'$], [$33, 43, 0', 4', 8'', 15''$], [$35, 42, 46, 50, 0', 2'$]]
        \item \{0=570, 1=31920, 2=485298, 3=1497276, 4=1268136\} \newline [[$0, 1, 2, 0', 0'', \infty$], [$0, 10, 22, 29, 39, 56$], [$8, 14, 21, 53, 0', 10''$], [$9, 25, 47, 50, 0', 4'$], \newline [$10, 20, 24, 0'', 4'', 13''$], [$12, 23, 45, 0', 2', 16''$], [$14, 16, 47, 56, 0'', 1''$], [$0', 1', 8', 7'', 9'', 12''$]]
        \item \{0=570, 1=31920, 2=452124, 3=1521672, 4=1276914\} \newline [[$0, 1, 2, 0', 0'', \infty$], [$0, 3, 11, 17, 39, 46$], [$10, 39, 44, 0', 4'', 18''$], [$11, 50, 55, 0', 2'', 17''$], \newline [$15, 21, 48, 0', 1', 6'$], [$19, 23, 25, 0', 1'', 13''$]]
        \item \{0=570, 1=31920, 2=464436, 3=1507536, 4=1278738\} \newline [[$0, 1, 2, 0', 0'', \infty$], [$0, 13, 23, 25, 39, 53$], [$3, 25, 40, 43, 49, 0'$], [$7, 23, 28, 30, 0'', 4''$], \newline [$12, 31, 44, 0', 1', 8'$], [$16, 38, 46, 0', 4', 8''$], [$19, 48, 0', 2', 3'', 17''$], [$24, 34, 41, 0'', 1'', 12''$]]
        \item \{0=570, 1=32376, 2=447678, 3=1521444, 4=1281132\} \newline [[$0, 1, 2, 0', 0'', \infty$], [$0, 9, 28, 32, 34, 47$], [$5, 29, 32, 41, 55, 0'$], [$8, 21, 36, 47, 0'', 4''$], \newline [$9, 23, 0', 4', 1'', 13''$], [$14, 22, 36, 0', 2', 14''$], [$15, 20, 52, 0', 1', 8'$], [$23, 39, 52, 0'', 2'', 16''$]]
        \item \{0=570, 1=35568, 2=426132, 3=1533984, 4=1286946\} \newline [[$0, 1, 2, 0', 0'', \infty$], [$0, 6, 19, 34, 38, 47$], [$3, 47, 52, 0'', 4'', 10''$], [$6, 14, 31, 36, 51, 0''$], \newline [$6, 35, 46, 0', 1', 8'$], [$13, 16, 23, 37, 0', 2'$], [$20, 21, 47, 0', 8'', 13''$], [$27, 36, 0', 4', 5'', 6''$]]
        \item \{0=570, 1=38760, 2=469908, 3=1511640, 4=1262322\} \newline [[$0, 1, 2, 0', 0'', \infty$], [$0, 4, 21, 26, 29, 37$], [$3, 21, 46, 52, 0', 8''$], [$4, 38, 39, 0'', 2'', 5''$], \newline [$12, 29, 50, 0', 2', 7''$], [$14, 34, 47, 53, 0'', 1''$], [$31, 36, 42, 49, 0', 4'$], [$0', 1', 8', 4'', 10'', 14''$]]
        \item \{0=570, 1=32376, 2=469908, 3=1498416, 4=1281930\} \newline [[$0, 1, 2, 0', 0'', \infty$], [$0, 4, 8, 16, 38, 39$], [$0, 6, 13, 41, 50, 55$], [$4, 7, 25, 31, 0', 4'$], \newline [$6, 35, 46, 0', 1', 8'$], [$16, 26, 0', 2', 10'', 11''$], [$22, 33, 38, 48, 0'', 4''$], [$24, 34, 47, 0', 2'', 5''$]]
        \item \{0=570, 1=32832, 2=466488, 3=1529880, 4=1253430\} \newline [[$0, 1, 2, 0', 0'', \infty$], [$0, 5, 7, 25, 27, 50$], [$4, 45, 53, 0'', 4'', 10''$], [$5, 23, 53, 0', 8'', 15''$], \newline [$6, 28, 44, 50, 54, 0''$], [$6, 35, 46, 0', 1', 8'$], [$10, 48, 0', 4', 7'', 9''$], [$19, 26, 47, 50, 0', 2'$]]
        \item \{0=570, 1=33288, 2=439470, 3=1520988, 4=1288884\} \newline [[$0, 1, 2, 0', 0'', \infty$], [$3, 47, 52, 0'', 4'', 10''$], [$4, 16, 21, 27, 37, 0''$], [$10, 44, 45, 0'', 1'', 12''$], \newline [$11, 25, 27, 50, 0', 4''$], [$15, 29, 41, 0', 2', 3''$], [$17, 32, 38, 42, 0', 4'$], [$0', 1', 8', 13'', 16'', 18''$]]
        \item \{0=570, 1=33288, 2=485982, 3=1520532, 4=1242828\} \newline [[$0, 1, 2, 0', 0'', \infty$], [$0, 4, 6, 18, 23, 52$], [$6, 35, 46, 0', 1', 8'$], [$11, 36, 55, 0'', 4'', 10''$], \newline [$12, 45, 48, 0', 1'', 18''$], [$16, 34, 42, 0', 4'', 15''$], [$26, 27, 49, 0', 4', 10'$], [$29, 36, 44, 0', 2', 5''$]]
        \item \{0=570, 1=33744, 2=447678, 3=1513692, 4=1287516\} \newline [[$0, 1, 2, 0', 0'', \infty$], [$3, 37, 41, 0', 2', 5'$], [$4, 7, 21, 32, 56, 0'$], [$8, 11, 20, 50, 0', 4''$], \newline [$15, 19, 38, 0', 1'', 2''$], [$26, 33, 43, 52, 0'', 4''$], [$29, 30, 52, 0', 4', 13'$], [$0', 1', 8', 13'', 16'', 18''$]]
        \item \{0=570, 1=34200, 2=461358, 3=1494996, 4=1292076\} \newline [[$0, 1, 2, 0', 0'', \infty$], [$3, 37, 41, 0', 2', 5'$], [$4, 9, 14, 40, 48, 0'$], [$5, 8, 16, 53, 0', 10''$], \newline [$5, 20, 30, 56, 0'', 4''$], [$25, 32, 55, 0', 4', 2''$], [$41, 48, 55, 0'', 1'', 8''$], [$0', 1', 8', 7'', 9'', 12''$]]
        \item \{0=570, 1=34200, 2=467856, 3=1504800, 4=1275774\} \newline [[$0, 1, 2, 0', 0'', \infty$], [$0, 8, 11, 16, 21, 49$], [$3, 26, 28, 0'', 2'', 16''$], [$4, 8, 19, 26, 55, 0'$], \newline [$7, 18, 45, 0', 1'', 8''$], [$11, 12, 17, 47, 0', 10''$], [$11, 18, 25, 0'', 4'', 13''$], [$39, 53, 0', 1', 6', 4''$]]
        \item \{0=570, 1=34656, 2=448020, 3=1507080, 4=1292874\} \newline [[$0, 1, 2, 0', 0'', \infty$], [$0, 10, 12, 29, 40, 47$], [$4, 7, 39, 0', 5'', 10''$], [$6, 38, 42, 45, 50, 0'$], \newline [$14, 30, 43, 0'', 1'', 8''$], [$15, 29, 0', 1', 4', 7''$], [$21, 27, 36, 47, 0', 8''$], [$25, 29, 48, 0'', 4'', 10''$]]
        \item \{0=570, 1=34656, 2=451782, 3=1525092, 4=1271100\} \newline [[$0, 1, 2, 0', 0'', \infty$], [$0, 4, 14, 31, 45, 53$], [$4, 6, 21, 24, 54, 0''$], [$12, 31, 44, 0', 1', 8'$], \newline [$17, 36, 0', 2', 3'', 8''$], [$18, 46, 54, 0', 5'', 13''$], [$22, 27, 33, 50, 0', 4'$], [$29, 30, 52, 0'', 4'', 13''$]]
        \item \{0=570, 1=34656, 2=452124, 3=1528056, 4=1267794\} \newline [[$0, 1, 2, 0', 0'', \infty$], [$0, 4, 26, 41, 48, 55$], [$0, 19, 22, 32, 33, 53$], [$10, 31, 43, 50, 0', 2'$], \newline [$13, 21, 55, 0', 1'', 9''$], [$15, 20, 52, 0', 1', 8'$], [$19, 35, 37, 46, 0'', 4''$], [$39, 47, 0', 4', 3'', 17''$]]
        \item \{0=570, 1=35112, 2=447336, 3=1506168, 4=1294014\} \newline [[$0, 1, 2, 0', 0'', \infty$], [$3, 8, 24, 30, 0', 5''$], [$3, 11, 29, 34, 55, 0''$], [$5, 16, 20, 0', 2'', 15''$], \newline [$15, 25, 39, 54, 0', 2'$], [$18, 41, 43, 0', 4', 10'$], [$26, 33, 40, 0'', 1'', 12''$], [$0', 1', 8', 7'', 9'', 12''$]]
        \item \{0=570, 1=35568, 2=450072, 3=1512552, 4=1284438\} \newline [[$0, 1, 2, 0', 0'', \infty$], [$0, 4, 26, 41, 48, 55$], [$3, 11, 30, 50, 0', 2'$], [$4, 14, 20, 25, 37, 0''$], \newline [$4, 17, 19, 0', 4'', 15''$], [$20, 41, 0', 4', 1'', 3''$], [$22, 42, 47, 0'', 4'', 13''$], [$39, 49, 56, 0', 1', 8'$]]
        \item \{0=570, 1=35568, 2=457938, 3=1515972, 4=1273152\} \newline [[$0, 1, 2, 0', 0'', \infty$], [$3, 26, 28, 0', 2', 16'$], [$5, 19, 24, 29, 56, 0''$], [$8, 29, 39, 45, 0', 4'$], \newline [$9, 21, 30, 38, 0', 5''$], [$10, 44, 45, 0'', 1'', 12''$], [$25, 43, 51, 0', 1'', 18''$], [$0', 1', 8', 4'', 10'', 14''$]]
        \item \{0=570, 1=36024, 2=450414, 3=1516884, 4=1279308\} \newline [[$0, 1, 2, 0', 0'', \infty$], [$0, 13, 21, 41, 44, 46$], [$3, 20, 31, 0', 2'', 18''$], [$6, 48, 0', 4', 1'', 13''$], \newline [$8, 12, 15, 39, 50, 0''$], [$14, 39, 44, 53, 0', 2'$], [$21, 29, 37, 0', 1', 8'$], [$29, 30, 52, 0'', 4'', 13''$]]
        \item \{0=570, 1=36480, 2=464436, 3=1513008, 4=1268706\} \newline [[$0, 1, 2, 0', 0'', \infty$], [$0, 7, 13, 41, 48, 56$], [$5, 36, 40, 0', 2', 5'$], [$6, 35, 46, 0', 1', 8'$], \newline [$6, 37, 44, 50, 0'', 1''$], [$7, 10, 26, 28, 52, 0'$], [$13, 15, 27, 0', 1'', 15''$], [$45, 50, 0', 4', 9'', 18''$]]
        \item \{0=570, 1=36936, 2=448362, 3=1537860, 4=1259472\} \newline [[$0, 1, 2, 0', 0'', \infty$], [$0, 8, 16, 18, 28, 35$], [$4, 15, 35, 0', 4', 13'$], [$5, 16, 26, 38, 39, 0'$], \newline [$10, 19, 35, 52, 0'', 1''$], [$12, 31, 44, 0', 1', 8'$], [$14, 29, 43, 0', 1'', 10''$], [$36, 48, 0', 2', 5'', 8''$]]
        \item \{0=570, 1=36936, 2=483588, 3=1515744, 4=1246362\} \newline [[$0, 1, 2, 0', 0'', \infty$], [$0, 4, 23, 26, 28, 36$], [$6, 22, 33, 39, 0', 4'$], [$8, 31, 33, 51, 0'', 1''$], \newline [$9, 19, 26, 0'', 4'', 13''$], [$15, 26, 30, 0', 2', 16''$], [$29, 36, 45, 48, 0', 10''$], [$0', 1', 8', 7'', 9'', 12''$]]
        \item \{0=570, 1=37392, 2=451440, 3=1518024, 4=1275774\} \newline [[$0, 1, 2, 0', 0'', \infty$], [$0, 9, 22, 28, 41, 56$], [$4, 13, 16, 27, 53, 0'$], [$6, 16, 31, 55, 0'', 1''$], \newline [$6, 35, 46, 0', 1', 8'$], [$8, 49, 54, 0'', 4'', 13''$], [$21, 24, 0', 4', 14'', 16''$], [$42, 47, 54, 0', 2', 11''$]]
        \item \{0=570, 1=37392, 2=469566, 3=1509588, 4=1266084\} \newline [[$0, 1, 2, 0', 0'', \infty$], [$3, 6, 24, 54, 0', 4'$], [$4, 8, 15, 29, 31, 0''$], [$4, 13, 21, 32, 0', 5''$], \newline [$5, 25, 27, 0'', 2'', 16''$], [$15, 25, 47, 0', 2', 1''$], [$40, 50, 54, 0'', 1'', 8''$], [$0', 1', 8', 4'', 10'', 14''$]]
        \item \{0=570, 1=37392, 2=465804, 3=1507080, 4=1272354\} \newline [[$0, 1, 2, 0', 0'', \infty$], [$0, 4, 13, 27, 44, 45$], [$4, 9, 42, 0', 4', 6''$], [$5, 25, 27, 0'', 2'', 16''$], \newline [$6, 35, 46, 0', 1', 8'$], [$12, 18, 21, 0', 1'', 8''$], [$27, 48, 55, 0', 5'', 15''$], [$31, 51, 56, 0', 2', 5'$]]
        \item \{0=570, 1=37848, 2=443574, 3=1525092, 4=1276116\} \newline [[$0, 1, 2, 0', 0'', \infty$], [$0, 5, 6, 28, 43, 53$], [$3, 19, 21, 52, 0', 2'$], [$6, 9, 26, 32, 47, 0''$], \newline [$8, 15, 22, 0'', 4'', 10''$], [$12, 31, 44, 0', 1', 8'$], [$16, 55, 0', 4', 1'', 13''$], [$25, 30, 42, 0', 2'', 18''$]]
        \item \{0=570, 1=38304, 2=466488, 3=1484280, 4=1293558\} \newline [[$0, 1, 2, 0', 0'', \infty$], [$0, 5, 14, 28, 34, 44$], [$6, 14, 22, 0'', 1'', 12''$], [$10, 19, 41, 0', 4', 5''$], \newline [$20, 23, 48, 0', 8'', 10''$], [$21, 44, 46, 0'', 4'', 13''$], [$22, 40, 43, 0', 2', 4''$], [$39, 49, 56, 0', 1', 8'$]]
        \item \{0=570, 1=38760, 2=461016, 3=1532160, 4=1250694\} \newline [[$0, 1, 2, 0', 0'', \infty$], [$0, 7, 21, 40, 50, 53$], [$4, 10, 15, 0', 2'', 10''$], [$5, 25, 27, 0'', 2'', 16''$], \newline [$6, 37, 0', 1', 6', 5''$], [$9, 32, 50, 56, 0', 1''$], [$12, 17, 29, 30, 38, 0'$], [$28, 32, 51, 0'', 4'', 13''$]]
        \item \{0=570, 1=38760, 2=469224, 3=1520304, 4=1254342\} \newline [[$0, 1, 2, 0', 0'', \infty$], [$0, 13, 23, 25, 39, 53$], [$5, 10, 13, 17, 37, 0'$], [$6, 35, 46, 0', 1', 8'$], \newline [$9, 24, 28, 0', 1'', 4''$], [$10, 29, 44, 53, 0'', 4''$], [$16, 27, 0', 4', 17'', 18''$], [$31, 36, 47, 0', 2', 16'$]]
        \item \{0=570, 1=38760, 2=473328, 3=1494768, 4=1275774\} \newline [[$0, 1, 2, 0', 0'', \infty$], [$0, 9, 38, 43, 49, 53$], [$4, 7, 25, 32, 52, 0'$], [$6, 37, 0', 1', 6', 5''$], \newline [$11, 18, 25, 0'', 4'', 13''$], [$12, 31, 44, 0'', 1'', 8''$], [$20, 29, 30, 55, 0', 2''$], [$42, 50, 53, 0', 1'', 15''$]]
        \item \{0=570, 1=39672, 2=432972, 3=1535808, 4=1274178\} \newline [[$0, 1, 2, 0', 0'', \infty$], [$0, 4, 23, 26, 28, 36$], [$4, 24, 29, 0', 2', 16'$], [$5, 48, 0', 4', 3'', 5''$], \newline [$6, 35, 46, 0', 1', 8'$], [$8, 27, 30, 45, 54, 0'$], [$9, 15, 28, 0', 4'', 10''$], [$9, 17, 21, 31, 0'', 1''$]]
        \item \{0=570, 1=39672, 2=466146, 3=1525092, 4=1251720\} \newline [[$0, 1, 2, 0', 0'', \infty$], [$0, 4, 21, 28, 47, 50$], [$5, 6, 25, 39, 50, 0'$], [$15, 29, 0', 1', 4', 7''$], \newline [$23, 39, 52, 0'', 2'', 16''$], [$24, 34, 41, 0'', 1'', 12''$], [$24, 40, 44, 0', 5'', 18''$], [$27, 33, 42, 45, 0', 10''$]]
        \item \{0=570, 1=39672, 2=465462, 3=1512324, 4=1265172\} \newline [[$0, 1, 2, 0', 0'', \infty$], [$0, 5, 25, 35, 37, 39$], [$0, 11, 22, 27, 40, 56$], [$9, 23, 28, 0', 4'', 15''$], \newline [$11, 48, 0', 2', 1'', 5''$], [$18, 25, 44, 52, 0'', 2''$], [$26, 41, 47, 50, 0', 4'$], [$39, 49, 56, 0', 1', 8'$]]
        \item \{0=570, 1=39672, 2=497610, 3=1489524, 4=1255824\} \newline [[$0, 1, 2, 0', 0'', \infty$], [$3, 23, 40, 0', 2', 9''$], [$5, 46, 51, 0'', 4'', 10''$], [$6, 50, 55, 0', 4', 13'$], \newline [$7, 14, 25, 29, 39, 0'$], [$10, 31, 40, 55, 0'', 1''$], [$20, 41, 42, 49, 0', 2''$], [$0', 1', 8', 13'', 16'', 18''$]]
        \item \{0=570, 1=41040, 2=475722, 3=1502748, 4=1263120\} \newline [[$0, 1, 2, 0', 0'', \infty$], [$0, 9, 29, 44, 49, 55$], [$3, 11, 30, 0', 2'', 17''$], [$5, 8, 9, 41, 51, 0'$], \newline [$11, 42, 46, 0'', 1'', 12''$], [$17, 26, 33, 56, 0', 1''$], [$18, 28, 35, 0'', 2'', 5''$], [$27, 32, 0', 1', 4', 13''$]]
        \item \{0=570, 1=41496, 2=460332, 3=1515744, 4=1265058\} \newline [[$0, 1, 2, 0', 0'', \infty$], [$0, 7, 15, 44, 50, 55$], [$5, 46, 51, 0'', 4'', 10''$], [$11, 21, 41, 44, 0', 10''$], \newline [$15, 25, 29, 47, 0'', 1''$], [$15, 28, 42, 48, 0', 4'$], [$24, 32, 33, 0', 2', 16''$], [$0', 1', 8', 7'', 9'', 12''$]]
        \item \{0=570, 1=41952, 2=499320, 3=1490664, 4=1250694\} \newline [[$0, 1, 2, 0', 0'', \infty$], [$0, 5, 8, 13, 21, 37$], [$0, 11, 28, 31, 33, 47$], [$5, 40, 49, 55, 0'', 4''$], \newline [$6, 35, 46, 0', 1', 8'$], [$16, 30, 0', 4', 17'', 18''$], [$20, 21, 44, 0', 1'', 4''$], [$25, 37, 48, 55, 0', 2'$]]
        \item \{0=570, 1=42408, 2=467172, 3=1503888, 4=1269162\} \newline [[$0, 1, 2, 0', 0'', \infty$], [$0, 6, 25, 29, 38, 40$], [$0, 7, 35, 41, 42, 55$], [$7, 23, 27, 47, 0', 4'$], \newline [$10, 16, 28, 0', 2'', 18''$], [$15, 39, 43, 51, 0'', 1''$], [$21, 29, 37, 0', 1', 8'$], [$31, 38, 0', 2', 11'', 17''$]]
        \item \{0=570, 1=43320, 2=491112, 3=1489296, 4=1258902\} \newline [[$0, 1, 2, 0', 0'', \infty$], [$0, 3, 13, 19, 42, 53$], [$5, 23, 26, 0', 2', 5''$], [$6, 35, 46, 0', 1', 8'$], \newline [$10, 47, 49, 0', 4'', 18''$], [$19, 39, 44, 0'', 4'', 10''$], [$19, 41, 50, 0', 4', 15''$], [$23, 28, 36, 0'', 1'', 8''$]]
        \item \{0=570, 1=43776, 2=468882, 3=1510044, 4=1259928\} \newline [[$0, 1, 2, 0', 0'', \infty$], [$0, 4, 8, 27, 43, 53$], [$3, 20, 50, 0', 2', 18''$], [$3, 26, 28, 0'', 2'', 16''$], \newline [$4, 25, 42, 55, 0', 4'$], [$5, 41, 46, 48, 0', 4''$], [$27, 42, 47, 48, 0'', 1''$], [$0', 1', 8', 2'', 11'', 15''$]]
        \item \{0=570, 1=43776, 2=478458, 3=1500924, 4=1259472\} \newline [[$0, 1, 2, 0', 0'', \infty$], [$3, 16, 34, 46, 0', 2'$], [$4, 15, 20, 21, 25, 0''$], [$9, 43, 47, 0'', 1'', 12''$], \newline [$23, 29, 31, 44, 0', 5''$], [$27, 48, 55, 0', 4', 8''$], [$32, 52, 54, 0'', 2'', 5''$], [$0', 1', 8', 2'', 11'', 15''$]]
        \item \{0=570, 1=44232, 2=471276, 3=1503888, 4=1263234\} \newline [[$0, 1, 2, 0', 0'', \infty$], [$0, 3, 14, 35, 44, 45$], [$6, 35, 46, 0', 1', 8'$], [$8, 45, 50, 0', 4', 8''$], \newline [$9, 26, 30, 0', 1'', 16''$], [$18, 48, 56, 0', 2'', 13''$], [$20, 27, 34, 0', 2', 5'$], [$30, 53, 55, 0'', 2'', 5''$]]
        \item \{0=570, 1=44232, 2=489744, 3=1492032, 4=1256622\} \newline [[$0, 1, 2, 0', 0'', \infty$], [$0, 8, 16, 19, 29, 33$], [$7, 33, 42, 0', 4', 5''$], [$11, 46, 48, 51, 0'', 2''$], \newline [$18, 25, 43, 55, 0', 2'$], [$22, 31, 46, 50, 0', 8''$], [$23, 28, 36, 0'', 1'', 8''$], [$0', 1', 8', 4'', 10'', 14''$]]
        \item \{0=570, 1=44688, 2=474696, 3=1512552, 4=1250694\} \newline [[$0, 1, 2, 0', 0'', \infty$], [$3, 6, 16, 49, 54, 0''$], [$3, 17, 28, 47, 53, 0'$], [$8, 27, 39, 0', 1'', 17''$], \newline [$21, 29, 37, 0', 1', 8'$], [$23, 43, 45, 0', 4', 13'$], [$29, 30, 52, 0'', 4'', 13''$], [$44, 56, 0', 2', 4'', 12''$]]
        \item \{0=570, 1=45144, 2=460674, 3=1513236, 4=1263576\} \newline [[$0, 1, 2, 0', 0'', \infty$], [$0, 8, 9, 23, 25, 31$], [$6, 17, 38, 47, 0', 4''$], [$8, 26, 32, 45, 0', 4'$], \newline [$9, 13, 35, 0', 2', 16'$], [$16, 34, 36, 0', 1'', 2''$], [$27, 44, 48, 55, 0'', 4''$], [$0', 1', 8', 13'', 16'', 18''$]]
        \item \{0=570, 1=45600, 2=489402, 3=1491348, 4=1256280\} \newline [[$0, 1, 2, 0', 0'', \infty$], [$0, 4, 12, 34, 38, 56$], [$5, 16, 33, 0'', 4'', 13''$], [$6, 25, 39, 54, 0', 4'$], \newline [$9, 34, 55, 0', 2', 4''$], [$17, 31, 33, 47, 0', 1''$], [$18, 39, 45, 55, 0'', 1''$], [$0', 1', 8', 13'', 16'', 18''$]]
        \item \{0=570, 1=46056, 2=457938, 3=1512324, 4=1266312\} \newline [[$0, 1, 2, 0', 0'', \infty$], [$0, 4, 12, 36, 46, 55$], [$11, 21, 39, 0', 4', 10''$], [$15, 20, 52, 0', 1', 8'$], \newline [$16, 42, 56, 0'', 2'', 16''$], [$26, 33, 40, 0'', 1'', 12''$], [$26, 42, 50, 0', 5'', 18''$], [$29, 30, 41, 0', 2', 3''$]]
        \item \{0=570, 1=46512, 2=473328, 3=1500240, 4=1262550\} \newline [[$0, 1, 2, 0', 0'', \infty$], [$0, 20, 22, 25, 33, 56$], [$4, 21, 24, 28, 30, 0'$], [$7, 11, 19, 0', 1'', 9''$], \newline [$22, 41, 51, 0', 2', 16'$], [$25, 34, 39, 52, 0'', 4''$], [$29, 54, 0', 4', 14'', 16''$], [$39, 49, 56, 0', 1', 8'$]]
        \item \{0=570, 1=46512, 2=482562, 3=1496820, 4=1256736\} \newline [[$0, 1, 2, 0', 0'', \infty$], [$3, 14, 30, 34, 0', 2'$], [$4, 7, 15, 0', 4'', 18''$], [$5, 19, 28, 32, 43, 0''$], \newline [$18, 41, 43, 0', 4', 10'$], [$39, 48, 51, 56, 0', 5''$], [$40, 50, 54, 0'', 1'', 8''$], [$0', 1', 8', 2'', 11'', 15''$]]
        \item \{0=570, 1=46968, 2=482562, 3=1501836, 4=1251264\} \newline [[$0, 1, 2, 0', 0'', \infty$], [$7, 17, 27, 33, 0', 1''$], [$11, 42, 46, 0'', 1'', 12''$], [$13, 15, 31, 36, 48, 0''$], \newline [$15, 34, 40, 0', 2'', 6''$], [$18, 28, 35, 0', 2', 5'$], [$36, 43, 51, 54, 0', 4'$], [$0', 1', 8', 13'', 16'', 18''$]]
        \item \{0=570, 1=47424, 2=454518, 3=1515060, 4=1265628\} \newline [[$0, 1, 2, 0', 0'', \infty$], [$0, 4, 26, 41, 48, 55$], [$0, 8, 18, 34, 38, 46$], [$4, 14, 37, 50, 0', 2'$], \newline [$12, 18, 43, 0', 2'', 12''$], [$22, 27, 30, 42, 0'', 2''$], [$34, 52, 0', 4', 1'', 13''$], [$39, 49, 56, 0', 1', 8'$]]
        \item \{0=570, 1=47424, 2=469566, 3=1506396, 4=1259244\} \newline [[$0, 1, 2, 0', 0'', \infty$], [$3, 16, 33, 53, 0', 2'$], [$3, 18, 34, 46, 54, 0''$], [$11, 29, 44, 0', 4', 16''$], \newline [$17, 19, 51, 0'', 1'', 8''$], [$20, 21, 45, 55, 0', 4''$], [$31, 36, 47, 0'', 2'', 16''$], [$0', 1', 8', 2'', 11'', 15''$]]
        \item \{0=570, 1=48792, 2=504792, 3=1475160, 4=1253886\} \newline [[$0, 1, 2, 0', 0'', \infty$], [$0, 5, 21, 27, 45, 52$], [$5, 15, 18, 0', 4', 14''$], [$12, 31, 44, 0', 1', 8'$], \newline [$17, 23, 54, 0', 4'', 12''$], [$19, 25, 42, 0', 1'', 17''$], [$26, 27, 49, 0'', 4'', 10''$], [$30, 38, 46, 0', 2', 16'$]]

        \item \{0=1140, 1=30096, 2=479826, 3=1511412, 4=1260726\} \newline [[$0, 1, 2, 0', 0'', \infty$], [$3, 6, 16, 41, 48, 0''$], [$8, 25, 31, 0', 4', 14''$], [$10, 20, 24, 0'', 4'', 13''$], \newline [$12, 23, 32, 44, 0', 2'$], [$21, 29, 37, 0'', 1'', 8''$], [$22, 30, 34, 51, 0', 5''$], [$0', 1', 8', 7'', 9'', 12''$]]
        \item \{0=1140, 1=31464, 2=448704, 3=1531704, 4=1270188\} \newline [[$0, 1, 2, 0', 0'', \infty$], [$3, 11, 20, 40, 46, 0''$], [$4, 47, 0', 4', 2'', 8''$], [$5, 6, 24, 48, 51, 0'$], \newline [$12, 31, 44, 0', 1', 8'$], [$13, 32, 37, 0', 1'', 15''$], [$14, 25, 42, 0', 2', 5'$], [$17, 19, 51, 0'', 1'', 8''$]]
        \item \{0=1140, 1=31920, 2=456228, 3=1524408, 4=1269504\} \newline [[$0, 1, 2, 0', 0'', \infty$], [$0, 5, 18, 32, 37, 49$], [$4, 11, 19, 47, 0', 4''$], [$8, 26, 46, 0', 5'', 12''$], \newline [$9, 34, 51, 55, 0', 2'$], [$23, 26, 27, 38, 0'', 2''$], [$25, 29, 48, 0', 4', 10'$], [$0', 1', 8', 2'', 11'', 15''$]]
        \item \{0=1140, 1=32832, 2=441522, 3=1527828, 4=1279878\} \newline [[$0, 1, 2, 0', 0'', \infty$], [$0, 6, 13, 41, 50, 55$], [$3, 8, 23, 26, 0', 4''$], [$5, 16, 33, 0'', 4'', 13''$], \newline [$7, 17, 25, 41, 47, 0'$], [$12, 29, 50, 0', 1'', 3''$], [$13, 32, 42, 0'', 1'', 8''$], [$27, 46, 0', 1', 4', 16''$]]
        \item \{0=1140, 1=34656, 2=457254, 3=1504572, 4=1285578\} \newline [[$0, 1, 2, 0', 0'', \infty$], [$0, 7, 10, 14, 33, 50$], [$4, 37, 39, 50, 55, 0'$], [$5, 39, 43, 52, 0'', 4''$], \newline [$13, 24, 44, 0'', 2'', 5''$], [$13, 25, 40, 0', 2', 1''$], [$15, 20, 52, 0', 1', 8'$], [$22, 43, 0', 4', 6'', 17''$]]
        \item \{0=1140, 1=34656, 2=465120, 3=1515744, 4=1266540\} \newline [[$0, 1, 2, 0', 0'', \infty$], [$5, 6, 9, 18, 24, 0''$], [$8, 15, 22, 0', 4', 10'$], [$9, 29, 37, 40, 0', 5''$], \newline [$12, 31, 33, 38, 0', 2'$], [$23, 39, 52, 0'', 2'', 16''$], [$34, 50, 52, 0', 1'', 8''$], [$0', 1', 8', 4'', 10'', 14''$]]
        \item \{0=1140, 1=34656, 2=466488, 3=1509360, 4=1271556\} \newline [[$0, 1, 2, 0', 0'', \infty$], [$0, 4, 13, 15, 35, 41$], [$0, 6, 34, 44, 49, 53$], [$3, 26, 40, 43, 0', 4'$], \newline [$10, 44, 0', 2', 5'', 13''$], [$15, 33, 36, 45, 0'', 2''$], [$18, 34, 50, 0', 4'', 8''$], [$39, 49, 56, 0', 1', 8'$]]
        \item \{0=1140, 1=34656, 2=495216, 3=1512096, 4=1240092\} \newline [[$0, 1, 2, 0', 0'', \infty$], [$3, 11, 40, 48, 0', 10''$], [$5, 9, 16, 20, 45, 0'$], [$6, 12, 37, 39, 0'', 1''$], \newline [$24, 39, 42, 0', 2', 16''$], [$26, 27, 49, 0', 4', 10'$], [$27, 31, 53, 0'', 4'', 13''$], [$0', 1', 8', 7'', 9'', 12''$]]
        \item \{0=1140, 1=36024, 2=452124, 3=1535352, 4=1258560\} \newline [[$0, 1, 2, 0', 0'', \infty$], [$0, 6, 19, 34, 38, 47$], [$0, 14, 20, 31, 40, 42$], [$3, 42, 47, 0', 2'', 10''$], \newline [$4, 34, 39, 55, 0'', 4''$], [$9, 33, 0', 4', 5'', 8''$], [$17, 20, 24, 51, 0', 2'$], [$39, 49, 56, 0', 1', 8'$]]
        \item \{0=1140, 1=36024, 2=460332, 3=1524408, 4=1261296\} \newline [[$0, 1, 2, 0', 0'', \infty$], [$0, 13, 23, 25, 39, 53$], [$5, 13, 20, 40, 56, 0'$], [$6, 35, 46, 0', 1', 8'$], \newline [$7, 38, 0', 4', 9'', 14''$], [$10, 14, 33, 0', 2', 16'$], [$20, 25, 44, 45, 0'', 4''$], [$41, 50, 53, 0', 1'', 12''$]]
        \item \{0=1140, 1=36480, 2=471276, 3=1502064, 4=1272240\} \newline [[$0, 1, 2, 0', 0'', \infty$], [$0, 11, 13, 23, 27, 52$], [$5, 12, 18, 44, 54, 0'$], [$6, 35, 46, 0', 1', 8'$], \newline [$16, 42, 56, 0', 2', 16'$], [$17, 18, 36, 53, 0'', 4''$], [$31, 52, 0', 4', 5'', 8''$], [$33, 38, 53, 0', 2'', 10''$]]
        \item \{0=1140, 1=36936, 2=455544, 3=1528512, 4=1261068\} \newline [[$0, 1, 2, 0', 0'', \infty$], [$0, 19, 21, 43, 53, 56$], [$9, 37, 40, 48, 0', 2'$], [$11, 18, 25, 0', 4', 13'$], \newline [$12, 27, 39, 50, 0', 5''$], [$15, 23, 55, 0', 2'', 10''$], [$30, 33, 39, 43, 0'', 4''$], [$0', 1', 8', 7'', 9'', 12''$]]
        \item \{0=1140, 1=37392, 2=481878, 3=1500468, 4=1262322\} \newline [[$0, 1, 2, 0', 0'', \infty$], [$0, 3, 9, 20, 45, 50$], [$8, 9, 13, 0', 8'', 10''$], [$15, 46, 48, 0', 1', 4'$], \newline [$19, 23, 51, 0', 1'', 14''$], [$26, 30, 38, 0', 4'', 5''$]]
        \item \{0=1140, 1=37848, 2=481536, 3=1510272, 4=1252404\} \newline [[$0, 1, 2, 0', 0'', \infty$], [$0, 10, 14, 17, 34, 36$], [$11, 23, 29, 36, 0', 5''$], [$17, 43, 54, 0'', 2'', 16''$], \newline [$20, 31, 47, 51, 0', 4'$], [$21, 44, 54, 0', 2', 1''$], [$23, 38, 39, 47, 0'', 1''$], [$0', 1', 8', 4'', 10'', 14''$]]
        \item \{0=1140, 1=38304, 2=478116, 3=1488840, 4=1276800\} \newline [[$0, 1, 2, 0', 0'', \infty$], [$0, 4, 17, 20, 22, 34$], [$9, 22, 50, 0', 2', 5'$], [$11, 15, 39, 0', 2'', 4''$], \newline [$14, 20, 52, 0', 1'', 13''$], [$20, 40, 42, 0'', 4'', 10''$], [$21, 29, 37, 0', 1', 8'$], [$36, 42, 45, 0', 4', 16''$]]
        \item \{0=1140, 1=38760, 2=466488, 3=1528056, 4=1248756\} \newline [[$0, 1, 2, 0', 0'', \infty$], [$0, 10, 13, 20, 24, 41$], [$6, 35, 46, 0', 1', 8'$], [$7, 13, 37, 50, 55, 0'$], \newline [$10, 14, 33, 0', 2', 16'$], [$20, 21, 47, 0', 1'', 13''$], [$28, 33, 47, 49, 0'', 2''$], [$48, 53, 0', 4', 9'', 18''$]]
        \item \{0=1140, 1=41040, 2=459648, 3=1518024, 4=1263348\} \newline [[$0, 1, 2, 0', 0'', \infty$], [$5, 10, 39, 45, 0', 4'$], [$7, 37, 40, 46, 51, 0''$], [$9, 25, 48, 0', 2'', 10''$], \newline [$12, 26, 43, 0', 2', 5'$], [$13, 30, 35, 49, 0', 5''$], [$23, 43, 45, 0'', 4'', 13''$], [$0', 1', 8', 7'', 9'', 12''$]]
        \item \{0=1140, 1=41040, 2=477090, 3=1509132, 4=1254798\} \newline [[$0, 1, 2, 0', 0'', \infty$], [$3, 47, 49, 0', 2', 4''$], [$4, 7, 20, 24, 51, 0'$], [$5, 20, 27, 43, 0'', 4''$], \newline [$7, 33, 47, 0'', 1'', 8''$], [$21, 27, 39, 42, 0', 1''$], [$25, 29, 48, 0', 4', 10'$], [$0', 1', 8', 13'', 16'', 18''$]]
        \item \{0=1140, 1=41040, 2=483588, 3=1500240, 4=1257192\} \newline [[$0, 1, 2, 0', 0'', \infty$], [$3, 6, 36, 49, 54, 0''$], [$4, 24, 29, 0'', 2'', 16''$], [$6, 21, 29, 40, 0', 4''$], \newline [$9, 16, 23, 45, 0', 2'$], [$20, 25, 55, 0', 5'', 12''$], [$26, 27, 49, 0', 4', 10'$], [$0', 1', 8', 2'', 11'', 15''$]]
        \item \{0=1140, 1=41496, 2=461700, 3=1527144, 4=1251720\} \newline [[$0, 1, 2, 0', 0'', \infty$], [$0, 8, 11, 16, 21, 49$], [$6, 10, 21, 30, 0'', 4''$], [$11, 20, 38, 0', 4', 14''$], \newline [$14, 34, 47, 48, 0', 5''$], [$15, 22, 29, 45, 0', 2'$], [$20, 31, 48, 0'', 1'', 12''$], [$0', 1', 8', 7'', 9'', 12''$]]
        \item \{0=1140, 1=41952, 2=481194, 3=1504116, 4=1254798\} \newline [[$0, 1, 2, 0', 0'', \infty$], [$0, 6, 15, 22, 25, 55$], [$9, 35, 47, 0', 2'', 9''$], [$10, 38, 54, 0', 4', 10'$], \newline [$11, 21, 41, 0', 5'', 18''$], [$15, 20, 52, 0', 1', 8'$], [$19, 29, 33, 0', 2', 5'$], [$36, 39, 50, 0', 1'', 15''$]]
        \item \{0=1140, 1=42864, 2=466146, 3=1503660, 4=1269390\} \newline [[$0, 1, 2, 0', 0'', \infty$], [$4, 34, 37, 43, 56, 0'$], [$7, 22, 53, 0', 4', 6''$], [$10, 38, 54, 0'', 4'', 10''$], \newline [$12, 17, 34, 50, 0'', 1''$], [$14, 16, 18, 28, 0', 1''$], [$30, 38, 46, 0', 2', 16'$], [$0', 1', 8', 13'', 16'', 18''$]]
        \item \{0=1140, 1=43776, 2=490428, 3=1479720, 4=1268136\} \newline [[$0, 1, 2, 0', 0'', \infty$], [$0, 5, 23, 25, 45, 55$], [$9, 13, 35, 0'', 2'', 16''$], [$10, 28, 33, 41, 0', 4'$], \newline [$11, 34, 40, 49, 0'', 1''$], [$15, 46, 54, 0', 2', 7''$], [$19, 22, 36, 43, 0', 8''$], [$0', 1', 8', 4'', 10'', 14''$]]
        \item \{0=1140, 1=43776, 2=491454, 3=1492716, 4=1254114\} \newline [[$0, 1, 2, 0', 0'', \infty$], [$0, 5, 9, 28, 39, 51$], [$5, 36, 40, 0'', 2'', 5''$], [$5, 48, 55, 0', 2', 12''$], \newline [$6, 35, 46, 0', 1', 8'$], [$13, 41, 54, 0', 1'', 2''$], [$16, 19, 40, 0', 4', 9''$], [$23, 43, 45, 0'', 4'', 13''$]]
        \item \{0=1140, 1=44688, 2=463752, 3=1505712, 4=1267908\} \newline [[$0, 1, 2, 0', 0'', \infty$], [$4, 6, 19, 36, 55, 0'$], [$5, 44, 45, 52, 0', 4''$], [$8, 49, 54, 0', 4', 13'$], \newline [$12, 26, 43, 0'', 2'', 5''$], [$13, 17, 38, 44, 0'', 1''$], [$18, 22, 46, 0', 2', 18''$], [$0', 1', 8', 2'', 11'', 15''$]]
        \item \{0=1140, 1=44688, 2=488376, 3=1504800, 4=1244196\} \newline [[$0, 1, 2, 0', 0'', \infty$], [$0, 3, 7, 29, 43, 53$], [$0, 10, 14, 17, 34, 36$], [$0, 20, 22, 25, 33, 56$], \newline [$5, 15, 54, 0', 4'', 8''$], [$9, 22, 28, 0', 1'', 17''$], [$13, 43, 52, 0', 2'', 13''$], [$39, 45, 50, 0', 1', 4'$]]
        \item \{0=1140, 1=45144, 2=469566, 3=1529652, 4=1237698\} \newline [[$0, 1, 2, 0', 0'', \infty$], [$3, 11, 23, 44, 0', 2'$], [$5, 25, 27, 0'', 2'', 16''$], [$7, 13, 30, 37, 45, 0''$], \newline [$10, 25, 39, 50, 0', 5''$], [$18, 28, 40, 0', 4'', 12''$], [$27, 31, 53, 0', 4', 13'$], [$0', 1', 8', 2'', 11'', 15''$]]
        \item \{0=1140, 1=45144, 2=496242, 3=1493172, 4=1247502\} \newline [[$0, 1, 2, 0', 0'', \infty$], [$0, 10, 21, 34, 38, 41$], [$5, 13, 0', 2', 13'', 16''$], [$6, 12, 17, 24, 39, 0''$], \newline [$6, 35, 46, 0', 1', 8'$], [$8, 31, 50, 52, 0', 4'$], [$20, 33, 37, 0', 4'', 8''$], [$40, 50, 54, 0'', 1'', 8''$]]
        \item \{0=1140, 1=49248, 2=503082, 3=1485876, 4=1243854\} \newline [[$0, 1, 2, 0', 0'', \infty$], [$4, 11, 29, 32, 41, 0'$], [$5, 40, 48, 56, 0', 4''$], [$13, 32, 42, 0'', 1'', 8''$], \newline [$15, 44, 55, 0', 2', 16'$], [$22, 26, 28, 0', 4', 16''$], [$23, 27, 30, 49, 0'', 2''$], [$0', 1', 8', 2'', 11'', 15''$]]

        \item \{0=1710, 1=32832, 2=454176, 3=1539912, 4=1254570\} \newline [[$0, 1, 2, 0', 0'', \infty$], [$3, 13, 25, 34, 50, 0''$], [$9, 25, 39, 0', 2', 9''$], [$9, 43, 47, 0'', 1'', 12''$], \newline [$10, 16, 43, 54, 0', 2''$], [$17, 24, 30, 53, 0', 4'$], [$27, 31, 53, 0'', 4'', 13''$], [$0', 1', 8', 13'', 16'', 18''$]]
        \item \{0=1710, 1=34200, 2=449046, 3=1533300, 4=1264944\} \newline [[$0, 1, 2, 0', 0'', \infty$], [$0, 19, 21, 43, 53, 56$], [$6, 35, 46, 0', 1', 8'$], [$7, 18, 33, 0', 1'', 9''$], \newline [$8, 13, 16, 20, 38, 0''$], [$16, 21, 0', 2', 14'', 17''$], [$21, 44, 46, 0'', 4'', 13''$], [$24, 31, 37, 52, 0', 4'$]]
        \item \{0=1710, 1=35112, 2=455886, 3=1529196, 4=1261296\} \newline [[$0, 1, 2, 0', 0'', \infty$], [$0, 3, 23, 25, 46, 56$], [$6, 16, 40, 46, 55, 0''$], [$9, 37, 56, 0'', 4'', 10''$], \newline [$12, 31, 44, 0', 1', 8'$], [$13, 16, 56, 0', 1'', 13''$], [$15, 23, 52, 54, 0', 4'$], [$27, 38, 0', 2', 5'', 8''$]]
        \item \{0=1710, 1=36936, 2=462726, 3=1523724, 4=1258104\} \newline [[$0, 1, 2, 0', 0'', \infty$], [$0, 5, 7, 25, 27, 50$], [$4, 40, 52, 56, 0'', 2''$], [$7, 10, 24, 36, 45, 0'$], \newline [$21, 29, 37, 0', 1', 8'$], [$23, 39, 52, 0', 2', 16'$], [$27, 47, 0', 4', 7'', 13''$], [$28, 43, 50, 0', 5'', 12''$]]
        \item \{0=1710, 1=37848, 2=465120, 3=1520760, 4=1257762\} \newline [[$0, 1, 2, 0', 0'', \infty$], [$3, 14, 37, 43, 52, 0''$], [$18, 44, 48, 51, 0', 2''$], [$20, 40, 42, 0', 4', 10'$], \newline [$21, 33, 38, 41, 0', 2'$], [$21, 44, 46, 0'', 4'', 13''$], [$26, 43, 56, 0', 1'', 9''$], [$0', 1', 8', 13'', 16'', 18''$]]
        \item \{0=1710, 1=45600, 2=485298, 3=1502292, 4=1248300\} \newline [[$0, 1, 2, 0', 0'', \infty$], [$0, 6, 25, 29, 38, 40$], [$5, 11, 12, 32, 0', 8''$], [$5, 25, 27, 0'', 2'', 16''$], \newline [$10, 26, 30, 0', 4', 5''$], [$15, 23, 41, 56, 0', 2'$], [$16, 23, 35, 40, 0'', 1''$], [$0', 1', 8', 4'', 10'', 14''$]]
        \item \{0=1710, 1=48792, 2=486666, 3=1503660, 4=1242372\} \newline [[$0, 1, 2, 0', 0'', \infty$], [$0, 4, 16, 18, 26, 53$], [$8, 18, 45, 50, 56, 0'$], [$11, 12, 34, 0'', 2'', 16''$], \newline [$15, 55, 0', 1', 4', 15''$], [$21, 33, 46, 54, 0', 8''$], [$27, 30, 44, 0', 4'', 17''$], [$41, 48, 55, 0'', 1'', 8''$]]
        \item \{0=1710, 1=57000, 2=466146, 3=1500468, 4=1257876\} \newline [[$0, 1, 2, 0', 0'', \infty$], [$0, 3, 10, 28, 36, 56$], [$11, 36, 52, 0', 5'', 10''$], [$12, 24, 50, 0', 1', 14'$], \newline [$13, 15, 33, 0', 1'', 8''$], [$22, 31, 37, 0', 2'', 6''$]]

        \item \{0=2280, 1=35568, 2=472644, 3=1515288, 4=1257420\} \newline [[$0, 1, 2, 0', 0'', \infty$], [$3, 33, 44, 53, 0', 4'$], [$3, 47, 52, 0'', 4'', 10''$], [$5, 6, 25, 29, 32, 0''$], \newline [$12, 26, 43, 0', 2', 5'$], [$13, 30, 38, 42, 0', 5''$], [$14, 50, 56, 0', 2'', 10''$], [$0', 1', 8', 7'', 9'', 12''$]]
    \end{enumerate}
    \end{example}

    \section{Acknowledgement}
    The authors express their thanks to Vedran Krčadinac for valuable remarks and verifying that our results are reproducible.

\end{document}